\documentclass[11pt,a4paper]{amsart}
\numberwithin{equation}{section}
\usepackage{rotating}
\usepackage{mathrsfs}
\usepackage{amscd,amssymb,amsopn,amsmath,amsthm,graphics,amsfonts,enumerate,verbatim,calc}
\usepackage[all]{xy}
\usepackage[utf8]{inputenc}
\usepackage{color}

\usepackage[a4paper,margin=1in]{geometry}
\usepackage[utf8]{inputenc}
\usepackage{amscd,amssymb,amsopn,amsmath,amsthm,graphics,amsfonts,enumerate,verbatim,calc}
\usepackage[T1]{fontenc}
\usepackage{lmodern}
\usepackage{enumitem}
\usepackage{hyperref}
\usepackage{microtype}\usepackage{hyperref}
\usepackage{geometry}
\usepackage{enumitem}
\usepackage{booktabs}
\usepackage{tikz-cd}
\hypersetup{colorlinks=true,linkcolor=blue,citecolor=blue,urlcolor=blue}
\newtheorem{theorem}{Theorem}[section]
\newtheorem{proposition}[theorem]{Proposition}
\newtheorem{corollary}[theorem]{Corollary}
\newtheorem{lemma}[theorem]{Lemma}
\newtheorem{fact}[theorem]{Fact}\newtheorem{claim}[theorem]{Claim}
\newtheorem{example}[theorem]{Example}\newtheorem{notation}[theorem]{Notation}

\newtheorem{remark}[theorem]{Remark}

\newtheorem{definition}[theorem]{Definition}
\newcommand{\Min}{\operatorname{Min}} 
\newcommand{\Cl}{\operatorname{Cl}}

\newcommand{\divisor}{\operatorname{div}}

\newcommand{\Spec}{\operatorname{Spec}}

\newcommand{\depth}{\operatorname{depth}}

\newcommand{\Frac}{\operatorname{Frac}}

\DeclareMathOperator{\edim}{edim}

\newcommand{\m}{\mathfrak m}
\newcommand{\n}{\mathfrak n}
\newcommand{\fm}{\mathfrak m}
\newcommand{\fn}{\mathfrak n}

\newcommand{\ZZ}{\mathbb Z}

\newcommand{\RR}{\mathbb R}

\newcommand{\FF}{\mathbb F}

\DeclareMathOperator{\ord}{ord}

\newcommand{\F}{\mathbb F}
\newcommand{\Aset}{\mathcal A}
\newcommand{\perf}{^{\infty}}

\newenvironment{PROOF}[2][\proofname.]
{\begin{proof}[#1]\renewcommand{\qedsymbol}{\textsquare$_{\rm #2}$}}
	{\end{proof}}
\begin{document}
	
\author[M. Asgharzadeh]{Mohsen Asgharzadeh}
\date{}

\address{M. Asgharzadeh}
\email{mohsenasgharzadeh@gmail.com}

	\title[ Regular models for F-coherent rings]
{Regular models for F-coherent rings }


 \begin{abstract}
Let $R$ be a Noetherian local ring of characteristic $p>0$ with
perfection $R^\infty$. We study   whether coherence of $R^\infty$ forces
$R^\infty\cong S^\infty$ for some regular Noetherian local ring $S$.
For a normal Noetherian local domain, coherence of $R^\infty$ forces
$R$ to be $p$-almost factorial, equivalently $R^\infty$ is a GCD
	domain; a finite regular extension $R\subseteq S$ with
$S^q\subseteq R$ gives the uniform bound $q\Cl(R)=0$.  In the
non-Noetherian direction we construct, for every prime $p$ and every
$d\geq2$, a perfect coherent local valuation domain of dimension $d$
that is not the perfection of any Noetherian ring. In the positive
direction we characterize coherence for affine semigroup local rings. Finally, we give an F-coherent ring without regular model.
\end{abstract}

	\maketitle

\section{Introduction}
\label{sec:intro}

The guiding question of this note is the following.
$$(Q):\quad
R^\infty \text{ coherent}
\quad\stackrel{?}\Longrightarrow\quad
R^\infty \cong S^\infty
\text{ for some regular Noetherian local ring } S 
$$
Here $R^\infty=\bigcup_{e\geq0}R^{1/p^e}$ is the perfection of a
Noetherian local ring $R$ of characteristic $p>0$. Following \cite{Shim}
$R$ is called $F$-coherent, if $R^\infty$ is coherent.
The unrestricted
version of $(Q)$ is posed  by  Shimomoto  motivated via the search for a
counterexample to the expectation that $F$-coherent should admit a ``regular model''. 
We show:
	\begin{theorem} 
	For every prime $p$ there is a countable normal Noetherian local
	domain $R$ with
	$\widehat R\simeq \F_p[[u^2,u^3,v]]$
	such that $R\perf$ is coherent and $R\perf\not\simeq S\perf$ for
	every regular Noetherian local ring $S$.
\end{theorem}

In particular, this shows that  F-coherent rings properly extend beyond purely inseparable extensions of regular rings. This author has occasionally wondered whether this phenomenon can occur even when the ring has an isolated singularity.
The paper proves two other results explicitly:

(A)
	For a normal Noetherian local domain $R$, coherence of $R^\infty$
	forces $R$ to be $p$-almost factorial: every divisor class is killed
	by a power of $p$.
	Equivalently, $R^\infty$ is a GCD domain. If moreover $R$ admits a
	finite extension $R\subseteq S$ with $S$ regular and
	$S^q\subseteq R$ for a single $q=p^e$, then $q\Cl(R)=0$. This is a necessary condition on any
	normal counterexample, not a sufficient one.
	
(B)
	For affine semigroup local rings, coherence of the perfection is
	equivalent to a simplicial condition on the cone together with
	finiteness of the $p$-primary part of the class group of the
	normalization, and in that case a finite regular extension with uniform
	Frobenius containment exists. Thus $(Q)$ has an
	affirmative answer in this case.

The remainder of the introduction explains how these results are
organized and what remains open.
The mechanism is the equivalence,
proved as Theorem~\ref{thm:gcd}, between $R^\infty$ being a GCD domain,
the existence for each height-one prime $P$ of a power $p^e$ with
$P^{(p^e)}$ principal, and the condition that $\Cl(R)$ be
$p$-primary torsion. The implication from coherence to GCD-ness uses
the fact that a coherent perfection of a Noetherian local domain is a
GCD domain, recalled from \cite[Section~2]{Simpson}. We also record that
when $R$ is excellent and reduced, $F$-coherence passes to the
normalization, so one may assume normality in the excellent case.

The uniform bound in (A) is Proposition~\ref{prop:uniform}. For the
four-generator monomial rings
$R_n=k[[u^n,u^{n-1}v,uv^{n-1},v^n]]$ with $n=p^e\geq4$ and their
normalizations $B_n$, it yields $n\Cl(B_n)=0$, and the surrounding
discussion verifies directly that $B_n$ is the normalization of the
nonnormal ring $R_n$. The same discussion warns that $(Q)$ asks
only for an isomorphism of perfections; one must not assume without
argument that a regular model produces a finite extension with a
uniform Frobenius containment. Almost factoriality is complemented by an exact criterion for coherence
itself, Proposition~\ref{prop:colon}. For a finitely generated ideal
$I$ and element $u$ of a Noetherian domain $R$, the chain of
ideals $J_e=(I^{[q]}:_R u^q)^{1/q}A$ inside $A=R^\infty$ is
eventually stationary precisely when $A$ is coherent. This makes
explicit that almost factoriality controls GCDs and divisor classes but
not the relations required for coherence, and it provides a concrete
test that later examples use.
Subsection 2.2 talks with several reductions for regular model.
For instance, suppose $R$ has a regular model. Then Proposition \ref{prop:catenary}
shows that
	$R$ is universally catenary.
 This is the motivation behind in Theorem 1.1.

The positive families in (B) are proved as follows.
Theorem~\ref{thm:toric} treats affine semigroup local rings;
Corollary~\ref{cor:binomial} extends this to arbitrary binomial ideals
at the origin over an algebraically closed field, and
Corollary~\ref{cor:completion} shows invariance of coherence under
completion in this class. Corollary~\ref{thm:hypersurface} classifies
primitive binomial hypersurfaces: the examples there include
$xy-z^n$, coherent exactly when $n$ is a power of $p$, and a
normal hypersurface whose normalization is regular but which is not
$F$-coherent in any characteristic.  Section 4 presents a system of computational examples that support $(Q)$. Also, we  use  lexicographic monomial valuations with value  excludes every Noetherian model. This
settles the non-Noetherian analogue but leaves $(Q)$ open.
Section 5 reconstruct some beautiful results of Shimomoto \cite{Shim} on F-singularity of F-coherent rings by  weakening assumptions like excellent, complete or even reducedness. In the final section we give the proof of Theorem 1.1. We also refer the reader to \cite{Ryo,Asg26,MPa,Elham} for some recent homological investigations of \(R^\infty\). I would like to thank my coauthors for many fruitful discussions concerning \(R^\infty\).

	\section{Almost factorial, and some reductions}
 	\subsection{Almost factorial}
 \begin{definition}
 	A normal Noetherian local domain $R$ is \emph{almost factorial}
 	if its divisor class group $\Cl(R)$ is torsion.
 	It is \emph{$p$-almost factorial} in the terminology of this note
 	if every class is killed by a power of $p$.
 	No uniform power is part of this definition.
 \end{definition}
 
 Equivalently, almost factoriality means that for each height-one
 prime $P$, some symbolic power
$ P^{(m)}=P^mR_P\cap R
$ is principal. The $p$-version requires $m=p^e$.
 For a local normal domain this is the usual $\mathbb Q$-factorial
 condition.
 We use normalized valuations $v_P:\Frac(R)^\times\to\ZZ$ for
 height-one primes of a normal Noetherian domain. They extend uniquely
 to $\Frac(R)^\infty$, by
$
 v_P(z)=p^{-e}v_P(z^{p^e}).
$ Normality gives the membership criterion
 $$
 z\in R\quad\Longleftrightarrow\quad
 v_P(z)\geq0\text{ for every height-one prime }P,
 \qquad z\in\Frac(R).
 $$
 
 \begin{theorem}\label{thm:gcd}
 Let $R$ be a normal Noetherian local domain of characteristic $p$,
 	and put $A=R^\infty$. The following are equivalent:
 \begin{enumerate}
\item $A$ is a GCD domain.
 	\item For every height-one prime $P\subset R$, there exist $e\geq0$
and $c\in R\setminus\{0\}$ such that
$$\divisor_R(c)=p^eP,
\qquad\text{equivalently}\qquad P^{(p^e)}=cR.
 $$
\item $\Cl(R)$ is $p$-primary torsion.
\end{enumerate}
 \end{theorem}
 
 \begin{proof}	\emph{$(1)\Rightarrow(2)$.}
Write $P=(a_1,\ldots,a_r)$ and let $d\in A$ be a greatest
common divisor of these generators. Such a thing exists, by GCD.
	Choose $q=p^e$ with $d^q\in R$.	For every height-one prime $Q\subset R$,
$A_{R\setminus Q}\cong(R_Q)^\infty	$
is a valuation domain: $R_Q$ is a DVR.
GCDs localize. One way to see this is to use
 $aA\cap bA=(ab/\gcd(a,b))A	$	and the fact that localization preserves finite intersections.
Consequently
$$
 	v_Q(d)=\min_i v_Q(a_i)=
 	\begin{cases}
 		1,&Q=P,\\
 		0,&Q\neq P.
 	\end{cases}
$$At $P$, the minimum is one because the $a_i$ generate $PR_P$;
at $Q\neq P$, at least one $a_i$ is outside $Q$.
It follows that
$\divisor_R(d^q)=qP$.
Then, the ideal characterized by these valuation inequalities
is exactly $d^qR=P^{(q)}$.
 	
\medskip\noindent
\emph{$(2)\Longleftrightarrow(3)$.}
 Height-one primes generate the divisor group, and each divisor
 	has finite support. Thus powers of $p$ killing the individual
 	prime classes kill every class, and conversely.
 	
 	\medskip\noindent
 	\emph{$(3)\Rightarrow(1)$.}
 	Take nonzero $a,b\in R$ and form the effective divisor
 $	D=\sum_{\operatorname{ht}P=1}\min\{v_P(a),v_P(b)\}\,P.$	Only finitely many terms are nonzero.
For some $q=p^e$, the divisor $qD$ is principal 
$qD=\divisor_R(c)$ with	$ c\in\Frac(R)^\times.$
Since $qD$ is effective, normality implies $c\in R$.
 Let $d:=c^{1/q}\in A$.
The valuations of $a^q/c$ and $b^q/c$ are nonnegative, so these
belong to $R$.
Hence $d$ divides $a,b$ in $A$.
If $h\in A$ is another common divisor, then
$v_P(h)\leq v_P(d)$ for every $P$.
Choose a power $Q$ of $p$, divisible by $q$, such that
$h^Q\in R$. The fraction
$(d/h)^Q=c^{Q/q}/h^Q\in\Frac(R)$	has nonnegative valuations everywhere, so it lies in $R$.
Therefore $d/h\in A$.
This proves that $d$ is a GCD of $a,b$ in $A$.
 	For arbitrary nonzero $a,b\in A$, first raise both to a common
 	$p$-power placing them in $R$, and then apply the inverse
 	Frobenius automorphism of $A$ to the resulting GCD.
 \end{proof}
 
 \begin{corollary}\label{cor:coh}
 If $R$ is a normal Noetherian local domain with coherent perfection,
then $R$ is $p$-almost factorial, and \ref{thm:gcd}(2) holds for every
height-one prime.
 \end{corollary}The almost-factorial conclusion in the excellent setting is also
 recorded in \cite[Remark~4.9]{Pat}.
 \begin{proof}
A coherent perfection of a Noetherian local domain is a GCD 	domain; see \cite[Corollary 3.7]{moh}. Apply Theorem~\ref{thm:gcd}.
 \end{proof}
 
 If the original ring $R$ is excellent and reduced but not normal,
 $F$-coherence makes it a domain and puts its finite normalization
 $B$ inside $R^\infty$, with $B^\infty=R^\infty$
 \cite[Section~3]{Shim}. That normalization is local and
 $p$-almost factorial. This is a statement about $\Cl(B)$, not
 an assertion that the nonnormal ring $R$ is almost factorial.

 \begin{proposition}\label{prop:uniform}
  Let $R\subseteq S$ be a finite extension of Noetherian local
 domains of characteristic $p$, with $R$ normal, $S$ regular,
 and $S^q\subseteq R$ for a fixed $q=p^e$. Then
$	q\Cl(R)=0. $
 \end{proposition}
 \begin{proof}
 The extension induces a bijection on spectra and preserves heights.
 For a height-one $P\subset R$, let $Q\subset S$ be the unique
 height-one prime above it.
Since $S$ is factorial, write $Q=tS$.
Let $r$ be the ramification index of the DVR extension
$R_P\subseteq S_Q$. Thus
$v_Q|_{\Frac(R)}=r\,v_P$.
Because $t^q\in R$, one has $r\mid q$, and
$	\divisor_R(t^q)=(q/r)P. $	There is no contribution at another height-one prime:
the corresponding prime of $S$ is different from $tS$.
 Consequently the order of $[P]$ divides $q/r$, hence divides $q$.
 The prime classes generate $\Cl(R)$.
 \end{proof}

 Almost factoriality controls GCDs, not all the relations required
 for coherence. The following criterion makes that distinction
 explicit.
 
 \begin{proposition}\label{prop:colon}
 Let $R$ be a Noetherian domain of characteristic $p$,
$A=R^\infty$, $I\subseteq R$ a finitely generated ideal,
and $u\in R$.
For $q=p^e$, define
$ 	C_e:=(I^{[q]}:_R u^q)$, $J_e:=C_e^{1/q}A. $
Here $I^{[q]}$ is generated by $q$-th powers of the elements
of $I$, and $C_e^{1/q}A$ is the ideal generated by all
$c^{1/q}$, $c\in C_e$. Then
$J_0\subseteq J_1\subseteq$ and
$	(IA:_A u)=\bigcup_{e\geq0}J_e.$
Each $J_e$ is finitely generated.
The ring $A$ is coherent if and only if, for every $I,u$,
this chain is eventually stationary.
 \end{proposition}
 
\begin{proof}
If $c\in C_e$, then $c^p\in C_{e+1}$, proving the inclusions.
Taking $q$-th roots of	$cu^q\in I^{[q]}$ shows $J_e\subseteq(IA:u)$.
Conversely, write
$ au=\sum_i f_i b_i$ where $ f_i\in I, a,b_i\in A.$
Choose $q=p^e$ placing $a,b_i$ in $R^{1/q}$.
Then $a^qu^q\in I^{[q]}$, hence $a\in C_e^{1/q}$.
Noetherianity of $R$ makes $C_e$, and therefore $J_e$,
finitely generated.
 If $A$ is coherent, the colon $(IA:u)$ is finitely generated;
its generators lie in a common $J_e$, so the chain stabilizes.
Conversely, stabilization makes all these colons finitely generated.
	For any finitely generated ideal and element of $A$, a common
Frobenius automorphism moves its generators and the element into
$R$. Hence all finitely generated ideal colons in $A$ are
finitely generated. This implies coherence: induct on the number of generators of 	an ideal and project its syzygies to the last coordinate. The image is the corresponding colon ideal, and the kernel is
the preceding syzygy module. Annihilators cause no additional condition.
 \end{proof}
 
 Simpson's factorial example \cite[Example~2.6]{Simpson}, with the obstruction
 made explicit:
 
 \begin{example}\label{ex:simpson}
 Let $p\equiv1\pmod5$ be prime and let 
 $	R:=\frac{\FF_p[[x_0,x_1,x_2,x_3,x_4]]}
{(x_0^5+x_1^5+x_2^5+x_3^5+x_4^5)}.$	This is a complete noetherian local UFD. Nevertheless its perfection
 	is a noncoherent GCD domain. 
 \end{example}
 
\begin{proof}
The ring has dimension four and an isolated singularity:
$p\neq5$, and the radical of its Jacobian ideal is the maximal ideal.
The hypersurface is normal by Serre's criterion.
Grothendieck's parafactoriality theorem for isolated local complete
	intersections of dimension at least four makes it factorial, as
in \cite[Example~2.6]{Simpson}.
Hence its perfection $A$ is a GCD domain by
Theorem~\ref{thm:gcd}.
Set 	$ I:=(x_1,x_2,x_3,x_4)R$,  $u;=x_0^4.$
The coefficient of $(x_0x_1x_2x_3x_4)^{p-1}$ in
$(\sum x_i^5)^{p-1}$ is
$	\frac{(p-1)!}{(((p-1)/5)!)^5}\neq0\pmod p.	$
Fedder's criterion \cite{Fedder} shows that $R$ is $F$-pure.
Thus $IA\cap R=I$.
Since $u\notin I$, the ideal
 $J=(IA:_A u) $
is proper.
For every $q=p^e=5m+1$, the defining relation gives	$ u^q=(-1)^{4m}x_0^4(x_1^5+x_2^5+x_3^5+x_4^5)^{4m}. $
In an expansion term whose four indices sum to $4m$, either
one index is at least $m+1$, or all four indices equal $m$.
The first type of term is in $I^{[q]}$; multiplication by $x_1$
runs the second type there as well.
 Consequently
$$(\ast):\quad	x_1u^q\in I^{[q]},
 \qquad x_1^{1/q}\in J
\quad\text{for every }q=p^e.$$
The ordinary order on $R$ is a valuation, since
 	$	\operatorname{gr}_{\mathfrak m}R
 	=\FF_p[x_0,\ldots,x_4]/(\textstyle\sum x_i^5) $	is a domain. Extend it to $A$ by taking $p$-power roots,
 	and call it $v$.
It is nonnegative on $A$ and strictly positive on every nonzero
nonunit. If $J=(g_1,\ldots,g_s)$, properness would imply
$	v(a)\geq\min_i v(g_i)>0$.
This contradicts  $(\ast)$, because $v(x_1^{1/q})=1/q\to0$.
Thus $J$ is not finitely generated, and $A$ is not coherent.
 \end{proof}
 
The following uses the fact that $F$-coherence descends along
 faithfully flat maps \cite[Section~3]{Shim}.
 
 \begin{example}\label{thm:candidate}
 	For every prime $p$, there exists a singular Noetherian local
 	UFD $R$ of characteristic $p$ with
 	$$
 (\ast):\quad	\quad\dim R=3,\qquad \edim R=4,\qquad R^\infty\text{ is coherent},	$$
 	and with completion
 $\widehat R\cong
T:=\frac{\FF_p[[X,Y,Z,W]]}{(XY-Z^p)}. $
In particular, $R$ is an abstract local hypersurface: its
completion is a codimension-one complete intersection.
Thus $R$ is factorial, hence almost factorial, and satisfies
all the input hypotheses for the embedding-dimension-four
version of $(Q)$.
 \end{example}
 
 \begin{proof}
The complete local domain $T$ has dimension three, depth three,
and embedding dimension four.
It is a hypersurface and every nonzero element of its prime field
is a nonzerodivisor.
Heitmann's completion theorem \cite{Heitmann} therefore supplies
a noetherian local UFD $R$ with $\widehat R\cong T$.
The canonical injection $R\to\widehat R\cong T$ shows that
$R$ has characteristic $p$.
Completion preserves dimension and embedding dimension, giving
the two equalities in $ (\ast)$.
Since $3<4$, $R$ is singular.
Embed $T$ in the complete regular ring
$S_0=\FF_p[[s,t,W]]$ by $X=s^p,\ Y=t^p,\ Z=st,\ W=W$.
This identifies $T$ with
$\FF_p[[s^p,t^p,st,W]]$, and
$S_0^p\subseteq T\subseteq S_0.$
 Hence $T^\infty=S_0^\infty$ is coherent.
 The completion map $R\to\widehat R$ is faithfully flat, so
 $F$-coherence descends to $R$, proving all claims.
 \end{proof}

 For excellent rings there is also a useful reduction:
 if $R$ were a counterexample, the normalization of $R_{\mathrm{red}}$ would
 have the same coherent perfection and would still be a counterexample.
 It would therefore be a normal $p$-almost-factorial Noetherian
 local domain by Corollary~\ref{cor:coh}.
 The almost-factorial conclusion is thus a necessary restriction on
 the search; it is not the missing proof of a Noetherian counterexample.
  \subsection{Reduction for regular model can be made }

  The following result concerns
 \emph{every} regular Noetherian model, including an initially
 noncomplete one. Coherence is not assumed in the theorem.
 
 \begin{lemma}\label{lem:parameter}
 Let $R$ be a complete Noetherian local domain of equal
 characteristic, with coefficient field $k$, dimension $d$,
and system of parameters $b_1,\ldots,b_d$.
Then the homomorphism
$k[[Y_1,\ldots,Y_d]]\to R$ via tha assignment $ Y_i\mapsto b_i	$ 	is injective and finite.
\end{lemma}
 \begin{proof}
The quotient $R/(b_1,\ldots,b_d)$ is finite-dimensional over $k$.
Lifting a basis and successively expanding in the $b_i$ shows	that $R$ is a finite module over the image of the displayed map;
the series converge because the parameter-ideal topology agrees
with the maximal-ideal topology.
	Its kernel is prime. Finiteness and $\dim R=d$ imply that the
quotient of the $d$-dimensional power series ring by this kernel
still has dimension $d$; hence the kernel is zero.
This is the usual parameter-subring form of the Cohen structure
theorem.
 \end{proof}
 
 \begin{proposition}
\label{thm:complete-model}	Let $R$ be a complete Noetherian local domain of characteristic
$p>0$, with finite residue field $k$, and let $A=R^\infty$.
 The following conditions are equivalent.
\begin{enumerate}
\item $A\cong S^\infty$ for a regular Noetherian local ring $S$.
 \item Inside $R$ there is a complete regular local subring
 	$	B\cong k[[Y_1,\ldots,Y_d]]$
 	such that $R$ is finite over $B$ and $R^Q\subseteq B$
	for some $Q=p^e$.
\item Inside $A$ there is a complete regular local ring $U$,
finite over $R$, such that $U^Q\subseteq R$ for some $Q=p^e$.
 	\end{enumerate}
\end{proposition}
 
 \begin{proof}
 If $d=0$, then $R=k$ and the assertions are immediate.
 Assume $d>0$.
 	
 \emph{$(1)\Rightarrow(2)$.}
 Transport $S$ into $A$ using the assumed isomorphism, so that
$S\subseteq A=S^\infty$.	Perfection preserves spectra and residue fields up to perfection.	Thus $\dim S=d$, and the residue field of $S$ is $k$: 	it embeds in the finite field $k$ and its perfection equals $k$.	The complete ring $R$ has its canonical coefficient field $k$.	This copy of $k$ is also contained in $S$. 	Indeed, if $a\in k$, some $a^{p^e}$ belongs to $S$. 	Since $k$ is finite, a further Frobenius power sends 	$a^{p^e}$ back to $a$.
 Choose regular parameters $t_1,\ldots,t_d$ of $S$, and choose	$q=p^a$ with $	b_i:=t_i^q\in R$ for all $1\leq i\leq d$.
 These elements form a system of parameters of $R$.
 In fact, the homeomorphisms$$\operatorname{Spec}R
	\ \longleftarrow\ \operatorname{Spec}A
 	\ \longrightarrow\ \operatorname{Spec}S$$
show that the only prime of $R$ containing all the $b_i$
is its maximal ideal. By Lemma~\ref{lem:parameter},
 $	B:=k[[b_1,\ldots,b_d]]\subseteq R $
is a complete regular local ring and $R$ is finite over $B$.
Here a formal series is evaluated in the complete ring $R$.
Let  $S_0:=B^{1/q}\subseteq A$.
It is a complete regular local ring with coefficient field $k$
and parameters $t_1,\ldots,t_d$.

 \begin{claim}\label{c26}
We claim that the original copy of $S$ in $A$ is contained in
this copy of $S_0$.
 \end{claim}

 \begin{PROOF}{\ref{c26}}	
An abstract
isomorphism of complete rings would not suffice.
First,
$	\bigcap_{N\geq1}(t_1,\ldots,t_d)^N A=0.$	To prove this, take the order valuation on	$\Frac(B)$ induced by total degree in the $b_i$;	in particular, $v(b_i)=1$.	Extend it to the finite field extension $\Frac(R)/\Frac(B)$,	and then to $\Frac(R)^\infty$.	Its valuation ring contains $R$, since $R$ is integral over $B$,	and hence it contains $A$.	Moreover $v(t_i)=1/q$.	Every nonzero element of $(t_1,\ldots,t_d)^N A$	therefore has value at least $N/q$, proving $	\bigcap_{N\geq1}(t_1,\ldots,t_d)^N A=0.$ 	By the Cohen structure theorem,	$\widehat S\cong k[[t_1,\ldots,t_d]]$.	There is a $k$-isomorphism	$	\beta:\widehat S\to S_0$	sending each parameter $t_i$ to the element with the same name.	For $s\in S$, choose polynomial truncations $s_N$ of its	expansion in $\widehat S$, of total degree less than $N$.	Then, in $A$,
	$ 	s-s_N\in(t_1,\ldots,t_d)^N A$ and $
	\beta(s)-s_N\in(t_1,\ldots,t_d)^N A. $
	By $\bigcap_{N\geq1}(t_1,\ldots,t_d)^N A=0	$, $s=\beta(s)$.
	Thus $S\subseteq S_0\subseteq A$, and so
	$A=S^\infty=S_0^\infty=B^\infty.$
 \end{PROOF} 

Choose finite $B$-module generators $r_1,\ldots,r_h$ of $R$.	By this, a common $Q=p^e$ puts all $r_j^Q$ in $B$.	Raising $r=\sum c_jr_j$, with $c_j\in B$, to its $Q$-th power	gives $r^Q\in B$. Hence $R^Q\subseteq B$, as required.

 \emph{$(2)\Rightarrow(3)$.}	Take $U=B^{1/Q}\subseteq A$.	Then $R\subseteq U$, $U^Q=B\subseteq R$, and $U$ is complete
and regular. Since $k$ is perfect, $B^{1/Q}$ is finite over $B$;	the same finite set of module generators also generates it over	the intermediate ring $R$.
 	
 	\emph{$(3)\Rightarrow(1)$.}	The inclusions $U^Q\subseteq R\subseteq U$ give	$R^\infty=U^\infty$.
 \end{proof}

 Here is a second setting in which any regular model can be made
 finite. It addresses polynomial quotients localized at closed
 points, without assuming their completions are the original rings.
 
 \begin{lemma}\label{lem:intermediate-finite}
 	Let $R$ be an excellent normal noetherian local domain of
 	characteristic $p$, with perfect residue field.
 	If
$R\subseteq V\subseteq R^\infty$
 	and $V$ is a regular noetherian local ring, then $V$ is finite
 	over $R$.
 \end{lemma}
 
 \begin{proof}	Write $K=\Frac(R)$, $L=\Frac(V)$, and $d=\dim R$.	The extension is integral and purely inseparable, and induces a	homeomorphism on spectra. Thus $\dim V=d$, the ideal	$\mathfrak m_RV$ is primary for $\mathfrak m_V$, and the residue	fields agree, since the residue field of $R$ is perfect.		The induced homomorphism $	\widehat R\longrightarrow\widehat V$	is finite. Indeed, $\widehat V/\mathfrak m_R\widehat V$ is	finite-dimensional over the common residue field. Lifting finitely	many generators and applying the complete module argument proves	finiteness. The two adic topologies agree, by the primary-ideal 	assertion. 	Excellence and normality make $\widehat R$ a normal domain; 	$\widehat V$ is a regular local domain.
 	Both have dimension $d$. The kernel of this finite homomorphism
 	is therefore zero. Set
 $H:=\Frac(\widehat R)$ and $ E:=\Frac(\widehat V).$	Then $E/H$ is finite, and $L$ embeds in $E$.	Excellence makes $R\to\widehat R$ regular, so $H/K$ is a	separable field extension. Consequently it is linearly disjoint	from the algebraic purely inseparable extension $L/K$.	Equivalently, $H\otimes_K L$ is a field, and it embeds in $E$.	It follows that	$[L:K]=[HL:H]\leq[E:H]<\infty.$	This also proves finiteness when $L/K$ was initially allowed	to be infinite, by considering its finite subextensions. 	Because $R$ is excellent, its normalization $N$ in the finite 	extension $L/K$ is finite over $R$. 	The integral intermediate ring $V$ is an $R$-submodule of $N$.	Noetherianity of $R$ now makes $V$ finite.
 \end{proof}
 
 \begin{proposition}\label{thm:finite-type-model}	Let $k$ be a finite field, and let $R$ be a normal local domain 	essentially of finite type over $k$, at a closed point. 	If $R^\infty\cong S^\infty$ for a regular Noetherian local ring $S$,	then there is a finite regular local extension 	$$
(\ast):\quad 	R\subseteq U\subseteq R^\infty,\qquad U^Q\subseteq R
 	\quad\text{for some }Q=p^e.	$$	Conversely, $(\ast)$ gives such a regular model.
 \end{proposition}
 
 \begin{proof}
Identify $S$ with a subring of $A=R^\infty$.	As in Theorem~\ref{thm:complete-model}, the finite field $k$	is contained in $S$. 	Write $R=C_{\mathfrak n}$, with $C$ a finite type $k$-domain, 	and choose finitely many $k$-algebra generators $c_i$ of $C$. 	For some $q=p^a$, all $c_i^q$ lie in $S$; hence $C^q\subseteq S$. 	If $h\in C\setminus\mathfrak n$, then $h^q\in S$ is outside	the maximal ideal of $A$, and therefore is a unit of $S$. 	It follows that $R^q\subseteq S$.	Thus
$	R\subseteq S^{1/q}\subseteq A.$	The ring $S^{1/q}$ is regular and noetherian.	The residue field of $R$ is finite, since the point is closed,	and $R$ is excellent.	Lemma~\ref{lem:intermediate-finite} proves that 	$U:=S^{1/q}$ is finite over $R$. 	A common $p$-power of finitely many $R$-module generators 	of $U$ belongs to $R$, giving the uniform desired $Q$. 	The converse follows by taking perfections.
 \end{proof}
 
  There is also an obstruction independent of divisor classes.
 
 \begin{proposition}\label{prop:catenary}	If a Noetherian local domain $R$ has	$R^\infty\cong S^\infty$ for a regular Noetherian local ring $S$,	then $R$ is catenary.
 \end{proposition}
 \begin{proof}	Perfection induces homeomorphisms of spectra, preserving the	order on primes, inclusions, and saturated chains.	The spectrum of $R$ therefore has the same chain-length	property as that of $S$.	A regular Noetherian local ring is Cohen--Macaulay and hence	catenary.
 \end{proof}
 
 Thus a noncatenary $F$-coherent Noetherian local domain would
 also answer $(Q)$ negatively. It maybe nice to produce  such a domain.

 \begin{proposition}\label{prop:catenary} 	Let $R$ and $S$ be Noetherian rings of characteristic $p>0$ with 	$R^\infty\cong S^\infty$. Then $R$ is universally catenary if and only if	$S$ is universally catenary. In particular, if $S$ is regular, then	$R$ is universally catenary.
 \end{proposition}
 
 \begin{proof}	There is a natural description	$R[T_1,\ldots,T_n]^\infty
\cong R^\infty[T_1^{1/p^\infty},\ldots,T_n^{1/p^\infty}].$	Indeed, each polynomial has finitely many coefficients and monomials,
so a common Frobenius power places it at a finite stage. The assumed	isomorphism therefore induces
$\bigl(R[T_1,\ldots,T_n]\bigr)^\infty
 	\cong
\bigl(S[T_1,\ldots,T_n]\bigr)^\infty.$	Perfection induces a homeomorphism on prime spectra, including for	nonreduced rings. Hence the spectra of the two polynomial rings are 	homeomorphic. Their prime posets have the same saturated chains, so	one polynomial ring is catenary if and only if the other is.	For a Noetherian ring, universal catenarity is equivalent to	catenarity of all its polynomial rings. Moreover, Noetherian	Cohen--Macaulay rings, in particular regular rings, are universally 	catenary.
 \end{proof}

 \section{A positive answer for affine semigroup rings}
 
 Let $k$ be a perfect field. Let $M$ be a finitely generated
 submonoid of a lattice, with $ M\cap(-M)=\{0\}$. Let $L$ be the group generated by $M$, and
 $ d:=\operatorname{rank}L. $ 
 Thus $M$ is a positive affine semigroup. 
 \begin{notation}
 Let $C:=\RR_{\geq0}M\subseteq L_{\RR}, 
 \overline M:=C\cap L, 
 G:=\ZZ[1/p], L_p:=L\otimes_{\ZZ}G,$ and define
 $ M_p:=\bigcup_{e\geq0}p^{-e}M\subseteq L_p.$ For a monoid $N$, write $k[N]$ for its monoid algebra, with
 monomials $X^u$, $u\in N$. Let
 $\mathfrak m:=(X^u:0\neq u\in M)$ and $ R:=k[M]_{\mathfrak m}.$
 The case $d=0$ is just $R=k$, so assume $d\geq1$.
 \end{notation}

 \begin{lemma}\label{lem:class}
 Let $\sigma_1,\ldots,\sigma_r\in L^*$ be the primitive integral	support forms of the facets of $C$, so that
 $	C=\{u\in L_{\RR}:\sigma_i(u)\geq0\text{ for all }i\}.$
 For $B:=k[\overline M]_{\overline{\mathfrak m}}$, where
 $\overline{\mathfrak m}$ is the monomial maximal ideal, there is	an exact sequence	$$	0\longrightarrow L
 	\xrightarrow{\ \sigma\ }\ZZ^r
 	\longrightarrow\Cl(B)\longrightarrow0,
 	\qquad
 	\sigma(u)=(\sigma_1(u),\ldots,\sigma_r(u)).$$
 \end{lemma}
 
 \begin{proof}	The global toric divisor sequence identifies	$\Cl(k[\overline M])$ with the cokernel of $\sigma$:	the coordinate divisors correspond to the facets, and the divisor of	$X^u$ has coefficients $\sigma_i(u)$. This is the usual toric 	divisor sequence; see \cite[Theorem~4.1.3]{CLS}. For completeness, localization at the monomial maximal ideal does	not change this class group. The localization map is surjective. 	Every global class has a representative that is a finitely generated
 $L$-graded divisorial fractional ideal $I$. If $I$ becomes
principal after localization, then
$	\dim_k I/\overline{\mathfrak m}I=1.$	Choose a monomial representing the nonzero homogeneous residue class.	Graded Nakayama, using any positive integral grading on $M$,	shows that this monomial generates $I$ globally.	Thus the localization map on class groups is injective as well.	Finally, $\sigma$ is injective because $C$ is full dimensional and	pointed.
 \end{proof}
 \begin{fact} If an excellent reduced Noetherian ring \(R\) has coherent
 perfection, its finite normalization \(B\) lies in \(R^\infty\), and
$B^\infty=R^\infty.$
 Here \(B\subseteq R^\infty\) means that every element of \(B\) has a
 \(p\)-power in \(R\), see \cite[3.8]{Shim}.\end{fact}
 \begin{theorem}\label{thm:toric}	With the notation above, the following are equivalent.	\begin{enumerate}	\item $R^\infty$ is coherent.		\item There is a $G$-basis $b_1,\ldots,b_d$ of $L_p$ such that	$	M_p=G_{\geq0}b_1+\cdots+G_{\geq0}b_d.$		\item Both of the following conditions hold:		$$ 		\overline M\subseteq M_p,		\qquad		\operatorname{coker}(\sigma:L\longrightarrow\ZZ^r)		\text{ is a finite }p\text{-group}.	 	$$		\item There is a finite extension of local domains $R\subseteq S$		such that $S$ is regular and $S^Q\subseteq R$ for some power		$Q$ of $p$.
 	\end{enumerate}	In particular, whenever $R^\infty$ is coherent, there  is a regular	local ring $S$ with $R^\infty=S^\infty$. 
 \end{theorem}
 
 \begin{proof}	First note that
$	R^\infty=k[M_p]_{\mathfrak m_p},$	where $\mathfrak m_p$ is the ideal of nonconstant monomials.	For $u\in L_p$, the following elementary membership test will be used:
 	$$
(+):\quad \quad 	X^u\in k[M_p]_{\mathfrak m_p}\quad\Longleftrightarrow\quad u\in M_p.	$$	Indeed, if $fX^u\in k[M_p]$ for some 	$f\in k[M_p]\setminus\mathfrak m_p$, the nonzero constant term 	of $f$ supplies a nonzero term of exponent $u$; no other term 	can cancel it.
 	 	\medskip\noindent
 	\emph{$(1)\Rightarrow(3)$.}
 	The normalization of $R$ is
$B=k[\overline M]_{\overline{\mathfrak m}}$, a finite extension.
The normalization has a unique point above the vertex, since every 	nonconstant monomial in $\overline M$ has a positive integral power 	in $M$. 	By the normalization theorem recalled above,
$B\subseteq R^\infty$ and $B^\infty=R^\infty$.	Applying (+) to its monomials gives $\overline M\subseteq M_p$.
 Moreover $B$ is excellent, normal and $F$-coherent. Thus
$\Cl(B)$ is $p$-primary torsion. 	By Lemma~\ref{lem:class} it is finitely generated, hence a finite
 	$p$-group.
 	
 	\medskip\noindent
 	\emph{$(3)\Rightarrow(2)$.}
 	The first condition of \ref{thm:toric}(3) gives
 $	M_p=C\cap L_p. $	To see this, multiply an element of $C\cap L_p$ by a suitable	power of $p$ to place it in $\overline M$, and then use	$\overline M\subseteq M_p$.
Tensoring $	0\longrightarrow L
\xrightarrow{\ \sigma\ }\ZZ^r
 \longrightarrow\Cl(B)\longrightarrow0,$ with $G$, the second condition of \ref{thm:toric}(3) yields
 an isomorphism
 $	\sigma_G:L_p\xrightarrow{\ \sim\ }G^r.$	Consequently $r=d$. Since the support inequalities define $C$,	this isomorphism takes $C\cap L_p$ onto $G_{\geq0}^d$.	The inverse images of the standard basis vectors give the desired	$b_i$. In particular $C$ is simplicial.
 	
 	\medskip\noindent
 \emph{$(2)\Rightarrow(4)$.}
 Choose monoid generators $m_1,\ldots,m_N$ for $M$.
 Their coordinates in the basis $b_i$ are nonnegative elements
 of $G$. Choose $q=p^e$ clearing all their denominators and set
 $	u_i=X^{b_i/q}$ and
 $S_0=k[u_1,\ldots,u_d]\subseteq k[M_p].$
 The $u_i$ are algebraically independent, and the choice of $q$	ensures $k[M]\subseteq S_0$.
 Each $b_i/q$ lies in $M_p$. Hence, for one common power $Q$
 of $p$, every $u_i^Q$ lies in $k[M]$. Therefore
 $S_0^Q\subseteq k[M], $
 and $S_0$ is finite over $k[M]$.	After localization at the vertices,
 $R\subseteq S:=k[u_1,\ldots,u_d]_{(u_1,\ldots,u_d)}$	is finite and $S^Q\subseteq R$. In fact, if a denominator has	nonzero constant term, its $Q$-th power lies in	$k[M]\setminus\mathfrak m$, so it is already invertible after	localizing $k[M]$.	This $S$ is regular.
 	
 \medskip\noindent
 \emph{$(4)\Rightarrow(1)$.}
 The inclusions $R\subseteq S\subseteq R^\infty$ imply
 $R^\infty=S^\infty$. The perfection of the regular ring $S$ is
 coherent.
 \end{proof}

 \begin{corollary}
 \label{cor:binomial}
 Let $k$ be algebraically closed of characteristic $p>0$.	Let $I\subseteq (X_1,\ldots,X_n)$ be a binomial ideal in
 	$k[X_1,\ldots,X_n]$, and let
 $R=\bigl(k[X_1,\ldots,X_n]/I\bigr)_{(X_1,\ldots,X_n)}.$	If $R^\infty$ is coherent, there is a regular local ring $S$
 	with $R^\infty\cong S^\infty$.

 \end{corollary}
 
 \begin{proof}
 Recall that $R_{\mathrm{red}}$ is a domain. Thus exactly one	minimal prime $P$ of $I$ survives localization at the origin,	and $R_{\mathrm{red}}\cong(k[\underline{X}]/P)_{(\underline{X})}$.	Over an algebraically closed field, minimal primes of binomial	ideals are binomial. Moreover, a prime binomial ideal, after deleting	the variables it contains and rescaling the remaining variables,	defines an affine semigroup domain \cite{ES}.	Since $P\subseteq(\underline{X})$, no surviving nonconstant monomial becomes 	a unit: otherwise a binomial with nonzero constant term would 	belong to $P$.	The semigroup is therefore positive, and the given point is its 	monomial vertex. Apply Theorem~\ref{thm:toric}; passing from $R$ 	to $R_{\mathrm{red}}$ does not change the perfection.
 \end{proof}
 
 \begin{corollary}\label{cor:completion}	For the semigroup local domains in Theorem~\ref{thm:toric},	$$
 	R^\infty\text{ is coherent}
 	\quad\Longleftrightarrow\quad
 	\widehat R^\infty\text{ is coherent}.	$$
 \end{corollary}
 
 \begin{proof}	One implication follows from faithfully flat descent along	$R\to\widehat R$.	For the other, Theorem~\ref{thm:toric} supplies a finite local	extension $R\subseteq S$, with $S$ regular and $S^Q\subseteq R$.	Completion is exact on finite $R$-modules, so it gives an injection
$	\widehat R\subseteq\widehat S.$	The $\mathfrak m_RS$-adic and maximal-ideal-adic topologies on	$S$ agree. Thus $\widehat S$ is complete regular.	By continuity of the $Q$-power map and closedness of	$\widehat R$ inside the finite $\widehat R$-module $\widehat S$,	one has $\widehat S^{\,Q}\subseteq\widehat R$.	Their perfections coincide.	For a possibly nonreduced binomial quotient, apply this argument to	$R_{\mathrm{red}}$. Excellence ensures that its completion is
reduced, and nilpotents do not affect perfection.
 \end{proof}
 
The {positive orthant} in \(\mathbb{R}^n\) is the set of all points whose coordinates are all nonnegative.

 \begin{corollary}\label{thm:hypersurface}
Let $k$ be a perfect field of characteristic $p>0$, and let	$$	R:=
 	\left(
 	\frac{k[X_1,\ldots,X_r,Y_1,\ldots,Y_s,W_1,\ldots,W_t]}
 	{(\,X_1^{a_1}\cdots X_r^{a_r}
 		-Y_1^{b_1}\cdots Y_s^{b_s}\,)}
 	\right)_{\mathfrak m},$$
 where $\mathfrak m$ is generated by all variables,
 $r,s\geq1$, all exponents are positive, and	$\gcd(a_1,\ldots,a_r,b_1,\ldots,b_s)=1.$
 Then $R^\infty$ is coherent exactly in the following cases:
 \begin{enumerate}
 \item $r=s=1$, with no further restriction;
\item $r=1,\ s\geq2$, and $a_1=p^e$ for some $e\geq0$;
\item $s=1,\ r\geq2$, and $b_1=p^e$ for some $e\geq0$.
 	\end{enumerate}
 In every affirmative case the perfection comes from a regular local
 ring. The same criterion holds for the completed hypersurface
 $\widehat R$.
 \end{corollary}
 
 \begin{proof}
 Condition $\gcd(a_1,\ldots,a_r,b_1,\ldots,b_s)=1$ makes the binomial primitive, so its quotient is a
positive affine semigroup domain. The variables $W_i$ contribute
independent free monoid factors.
Suppose first $r,s\geq2$. 
Let \(V=\mathbb{R}^{r+s}/\mathbb{R}v\), where
\[
v=(a_1,\dots,a_r,-b_1,\dots,-b_s),\qquad a_i,b_j>0.
\]
Let \(\bar e_i\) (for \(i=1,\dots,r\)) and \(\bar f_j\) (for \(j=1,\dots,s\)) be the images of the standard coordinate vectors. The cone \(C\) is generated by these \(r+s\) vectors.

\begin{claim}\label{pd}
\(\bar e_i\) are distinct and  extremal, i.e.\ \(\bar e_i\) is not a nonnegative combination of the others.
\end{claim}

\begin{PROOF}{\ref{pd}}
Suppose
$\bar e_i=\sum_{k\ne i}\lambda_k\bar e_k+\sum_{j}\mu_j\bar f_j$ where $\lambda_k,\mu_j\ge 0.$	Lifting to \(\mathbb{R}^{r+s}\), this means
$e_i-(\sum_{k\ne i}\lambda_k e_k+\sum_j\mu_j f_j)=c\,v$
for some \(c\in\mathbb{R}\). Comparing coordinates:
\begin{itemize}
\item {\(i\)-th \(X\)-coordinate:} \(1=c\,a_i\), so \(c=1/a_i>0\).
\item {\(k\)-th \(X\)-coordinate (\(k\ne i\)):} \(-\lambda_k=c\,a_k>0\), so \(\lambda_k<0\) --- contradiction to \(\lambda_k\ge0\).
\end{itemize}
So \(\bar e_i\) is extremal. The \(Y\)-case is symmetric (use a coordinate where \(v\) is negative). Thus all \(r+s\) rays are extremal. 
Two coordinate vectors could map to the same ray only if one were a positive multiple of the other modulo \(\mathbb{R}v\). But \(\bar e_i\) and \(\bar e_k\) have the same sign pattern in the \(X\)-coordinates and zero in the \(Y\)-coordinates; a positive multiple would force matching sign patterns, which fails for \(\bar e_i\) vs.\ \(\bar f_j\), and for \(\bar e_i\) vs.\ \(\bar e_k\) with \(i\ne k\) the \(i\)-th coordinate comparison gives a contradiction. So the rays are distinct.
\end{PROOF}

The cone has dimension $r+s-1$ but $r+s$ extremal rays,	so it is not simplicial. Adding the $W_i$ does not change this
failure. Theorem~\ref{thm:toric} excludes coherence.
Now suppose $r=1,\ s\geq2$, put $c=a_1$, and omit the	independent $W_i$ temporarily. 
Now \(r=1\), so \(v=(c,-b_1,\dots,-b_s)\) with \(c=a_1>0\), \(b_j>0\). Let $w:=(\tfrac{b_1}{c},\dots,\tfrac{b_s}{c})$. The lattice is
$L=\mathbb{Z}^s+\mathbb{Z}w\subseteq \mathbb{Q}^s,$ realized in \(\mathbb{R}^s\) (the \(X_1\)-coordinate has been eliminated using the relation).  The monoid is generated by
$\{e_i, w\}_i.$

\begin{claim}\label{po}	\(\mathbb{R}_{\ge0}M=\mathbb{R}_{\ge0}^s\), i.e., the positive orthant.
\end{claim}

\begin{PROOF}{\ref{po}}	Clearly \(M\subseteq \mathbb{R}_{\ge0}^s\), since all generators have nonnegative coordinates, so \(\mathbb{R}_{\ge0}M\subseteq\mathbb{R}_{\ge0}^s\). Conversely, each standard basis vector \(e_i\in M\), so \(\mathbb{R}_{\ge0}^s=\mathbb{R}_{\ge0}\{e_1,\dots,e_s\}\subseteq\mathbb{R}_{\ge0}M\). 
\end{PROOF}

 Since the cone is the orthant,
$M\cap\mathbb{R}_{\ge0}e_i = M\cap(\text{the }i\text{-th coordinate ray}).$ An element of \(M\) on this ray has the form \(t e_i\) with \(t\ge0\). But elements of \(M\) are nonnegative integer combinations of \(e_k\)'s and \(w\); the \(j\)-th coordinate (\(j\ne i\)) is \(\lambda_j + \lambda_w b_j/c\). For this to vanish with \(\lambda_j,\lambda_w\ge0\) and \(b_j/c>0\), we need \(\lambda_j=\lambda_w=0\). Then the element is \(\lambda_i e_i\) with \(\lambda_i\in\mathbb{Z}_{\ge0}\), i.e.\ \(t\in\mathbb{N}\). So \(M\cap\mathbb{R}_{\ge0}e_i=\mathbb{N} e_i\). The same argument over \(G=\mathbb{Z}[1/p]\) gives \(M_p\cap\mathbb{R}_{\ge0}e_i=G_{\ge0}e_i\).

This is exactly what the proof uses: the free \(G\)-generators of \(M_p\) must lie on the extremal rays, and on the \(i\)-th ray the \(G\)-span is \(G e_i\), forcing \(L_p=G^s\) and hence \(b_i/c\in G\).
After eliminating the single relation the cone lives in \(\mathbb{R}^s\), and the claim is
$\mathbb{R}_{\ge0}M=\mathbb{R}_{\ge0}^s,$
i.e.\ the cone spanned by the monoid \(M\) is the whole positive orthant. This matters because:
\begin{itemize}
\item the extremal rays of \(\mathbb{R}_{\ge0}^s\) are precisely the coordinate rays \(\mathbb{R}_{\ge0}e_i\);
\item the monoid \(M\) meets the \(i\)-th coordinate ray only in \(\mathbb{N}e_i\) (and \(M_p\) only in \(G_{\ge0}e_i\));
\item so if \(M_p\) is a free \(G\)-monoid (condition (2) of Theorem~\ref{thm:toric}), its generators must lie on these coordinate rays, which forces \(L_p=G^s\) and hence \(b_i/c\in G\).
\end{itemize}
That is exactly the step   to conclude \(c=p^e\).
Conversely, if $c=p^e$, the regular local subring
 $	A=k[Y_1,\ldots,Y_s,W_1,\ldots,W_t]_{(Y,W)} $	has a finite purely inseparable extension $A\subseteq R$, because
 $X_1^c=\prod Y_j^{b_j}$. Therefore	$R^\infty=A^\infty$, which is coherent.	The case $s=1,\ r\geq2$ is symmetric.
 	
Finally, if $r=s=1$, put $a=a_1$, $b=b_1$.
The parametrization $X=u^b,\ Y=u^a$ identifies $R$ with
$k[u^a,u^b,W]_{(u^a,u^b,W)}.$
 Since $\gcd(a,b)=1$, the numerical semigroup generated by $a,b$
contains every sufficiently large integer. Choose a power	$q=p^e$ in that semigroup. Then	$ R\subseteq S=k[u,W]_{(u,W)}$ and $S^q\subseteq R. $ 
This gives the required regular model.
The assertion about completion follows from
Corollary~\ref{cor:completion}.
 \end{proof}
 
 \begin{example}\label{ex:An}
 	For $	R_n:=\bigl(k[x,y,z]/(xy-z^n)\bigr)_{(x,y,z)},$	one has
 $$	R_n^\infty\text{ coherent}
 \quad\Longleftrightarrow\quad n=p^e\quad(e\geq0).	$$	For $n=p^e$, a regular extension is
 $	S=k[u,v]_{(u,v)}	$ with	$x=u^n, y=v^n,z=uv,	$	and $S^n\subseteq R_n$.
 For $n>1$ these are singular hypersurfaces with regular perfections.	When $n$ is not a power of $p$, the failure is coherence itself;
 such a ring is not a counterexample to $(Q)$.
 \end{example}
 
 \begin{example}
The hypersurface
$	\bigl(k[x,y,z,w]/(xy-zw)\bigr)_{(x,y,z,w)} $
 is not $F$-coherent in any characteristic. Its cone has four
 extremal rays and dimension three, as in the first case of
Theorem~\ref{thm:hypersurface}.
 \end{example}
 
 \begin{example}\label{ex:normalization}
 Consider
 $	R=\bigl(k[x,y,z]/(z^6-x^2y^3)\bigr)_{(x,y,z)}.$	It is not $F$-coherent for any prime characteristic, by ~\ref{thm:hypersurface}.
 Its normalization is nevertheless the regular ring
 $S=k[u,v]_{(u,v)}$, with
 $x=u^3$, $y=v^2$ and $z=uv. $ 	Indeed $u,v$ are integral and	$v=z^3/(xy)$, $u=z/v$ belong to the fraction field of $R$.
 This normalization is not purely inseparable in the required
 ring-theoretic sense: a pure $u$-power belongs to $R$ only when
 its exponent is divisible by $3$, and a pure $v$-power only
 when its exponent is divisible by $2$.
 \end{example}
 
 \section{some examples}

  \subsection {A perfect coherent ring with no regular model}

 The counterexample below is a \emph{non-Noetherian} local ring of
 characteristic $p>0$ and dimension two. Its perfection is coherent but
 is not the perfection of any Noetherian ring, hence of no regular
 Noetherian ring. This does not answer the corresponding question for
 Noetherian local rings, posed in \cite{Shim}.  
 
 \begin{lemma}\label{lem:obstruction}
 	Let $A$ be a Noetherian ring of characteristic $p>0$.
 	If $A^\infty$ is a valuation domain, then $\dim A\leq 1$.
 \end{lemma}
 
 \begin{proof}
 Replacing $A$ by $A_{\mathrm{red}}$ does not change its dimension
or its perfection. We may therefore identify $A$ with a subring of
the domain $B:=A^\infty$.
 The extension $A\subseteq B$ is integral and induces a homeomorphism
$\Spec B\longrightarrow\Spec A$. Indeed, every element of $B$	has a $p$-power in $A$, so the unique prime above
$\mathfrak q\in\Spec A$ consists of those $b\in B$ for which	$b^{p^e}\in\mathfrak q$ for some $e\geq0$.
Consequently, $A$ is a local domain, say with maximal ideal
$\mathfrak n$, and
$\dim A=\dim B$ and $ \sqrt{\mathfrak nB}=\mathfrak m_B.$	If $\mathfrak n=0$, then $\dim A=0$.
Otherwise, write $\mathfrak n=(a_1,\ldots,a_r)$.
In a valuation domain, a finitely generated ideal is generated by one
of its generators. Thus $\mathfrak nB=a_iB$ for some $i$, and	$\sqrt{a_iB}=\mathfrak m_B.$	Since $A\subseteq B$ is integral, contraction gives	$\sqrt{a_iA}=\sqrt{a_iB}\cap A=\mathfrak n.$	The principal ideal theorem yields
$\dim A=\operatorname{ht}\mathfrak n\leq1$.
 \end{proof}
 
 \begin{example}\label{thm:example}
 For every prime $p$, there exists a local domain $V$ of
characteristic $p$ satisfying
$$\dim V=2,\qquad V^\infty=V,\qquad V\text{ is coherent},$$
 such that $V\not\cong A^\infty$ for every noetherian ring $A$.
 \end{example}
 
 \begin{proof}
 Let $k$ be a perfect field of characteristic $p$, for example
 $k=\mathbb F_p$. Inside a fixed algebraic closure of $k(x,y)$, define
$K_e:=k\bigl(x^{1/p^e},y^{1/p^e}\bigr)$ and $K:=\bigcup_{e\geq0}K_e,	$ 
where $x,y$ are algebraically independent.
The field $K$ is perfect.	Set $G=\ZZ[1/p]$ and equip $\Gamma=G\oplus G$ with the	lexicographic order. Define a valuation on $K_e$ by setting, for a	nonzero polynomial, 	$$ 	v\!\left(\sum_{i,j}c_{ij}x^{i/p^e}y^{j/p^e}\right)	=	\min_{\mathrm{lex}}
 	\left\{\left(\frac{i}{p^e},\frac{j}{p^e}\right):c_{ij}\neq0\right\},$$
 and extending to rational functions by
$v(f/g)=v(f)-v(g)$.
These valuations are compatible as $e$ increases. They therefore
define a valuation $v:K^\times\longrightarrow\Gamma$ with value
group $\Gamma$. Let
$$V:=\{0\}\cup\{a\in K^\times:v(a)\geq(0,0)\}.$$	This is a local valuation domain with fraction field $K$ and maximal	ideal $\mathfrak m=\{0\}\cup\{a\in K^\times:v(a)>(0,0)\}.	$	We verify its dimension directly. Define 	$$	P:=\{0\}\cup
 	\{a\in V\setminus\{0\}:\text{the first coordinate of }v(a)
 	\text{ is positive}\}.	$$
It is a prime ideal, and	$	(0)\subsetneq P\subsetneq\mathfrak m, $	since $x\in P$ and $y\in\mathfrak m\setminus P$.	If $a\in\mathfrak m\setminus\{0\}$ has value	$v(a)=(\alpha,\beta)$, then	$$	\sqrt{aV}=	\begin{cases}		P,&\alpha>0,\\
 		\mathfrak m,&\alpha=0\ \text{and}\ \beta>0.
\end{cases}	$$	For the first case, a sufficiently large positive integer multiple of	any positive first coordinate exceeds $\alpha$. For the second case,	the same assertion holds for positive second coordinates, while any 	positive first coordinate already exceeds $(0,\beta)$. 	Every nonzero prime therefore contains $P$; if it strictly contains 	$P$, it contains an element of $\mathfrak m\setminus P$ and hence 	equals $\mathfrak m$. Thus these are all the prime ideals, and 	$\dim V=2$.	Since $K$ is perfect, every $0\neq a\in V$ has a unique $p$-th root 	in $K$, and 	$	v(a^{1/p})=\frac{1}{p}v(a)\geq(0,0).$	Hence $a^{1/p}\in V$, so Frobenius on $V$ is bijective and	$V^\infty=V$.	Every nonzero finitely generated ideal of $V$ is principal and,	because $V$ is a domain, is free of rank one. Thus every finitely 	generated ideal is finitely presented, proving coherence. 	Finally, if $V\cong A^\infty$ for a Noetherian ring $A$, 	Lemma~\ref{lem:obstruction} would imply 	$\dim V=\dim A\leq1$, contradicting $\dim V=2$.
 \end{proof}
 
 \begin{remark}	The ring $V$ is necessarily non-Noetherian. More explicitly, its	maximal ideal cannot be principal: if $\mathfrak m=aV$ with	$a\neq0$, then $a^{1/p}\in\mathfrak m$, but	$	v(a^{1/p})=\frac{1}{p}v(a)<v(a), $	so $a^{1/p}\notin aV$. Since every finitely generated ideal of $V$	is principal, $\mathfrak m$ is not finitely generated.
 \end{remark}
 
 \begin{remark} 	The same construction works in every dimension $d\geq2$:	replace $k(x,y)$ by $k(x_1,\ldots,x_d)$, take its perfection,	and use the lexicographic monomial valuation with value group	$\ZZ[1/p]^d$.	Its valuation ring is perfect and coherent of dimension $d$,	and Lemma~\ref{lem:obstruction} excludes every Noetherian model.
 \end{remark}
 
 \subsection{Noetherian example}

 The binomial conclusion is affirmative under the hypotheses above;
 there is no counterexample there. For arbitrary complete
 intersections or arbitrary non-binomial hypersurfaces, this note
 does not establish a general affirmative answer or a counterexample
 retaining coherence. The unrestricted question $(Q)$.

 \begin{example}\label{prop:tower}	Let $A=k[X_1,\ldots,X_d]_{(X)}$, where $k$ is perfect.	Let $q_i=p^{e_i}$ and choose	$$f_i\in k[X_1,\ldots,X_d,Z_1,\ldots,Z_{i-1}], 	\qquad f_i(0)=0. $$
Localize at the origin the ring
$R:= 	A[Z_1,\ldots,Z_c]/ 	(Z_1^{q_1}-f_1,\ldots,Z_c^{q_c}-f_c).$ 	Then $R$ is a local complete intersection, finite free over $A$, 	and
 $R^\infty=A^\infty.$	In particular its perfection is coherent and has a regular 	Noetherian model. This also holds after completion.
 \end{example}
 
 \begin{proof}	Successive division by the monic equations gives an $A$-basis	$	Z_1^{j_1}\cdots Z_c^{j_c},$ for $ 0\leq j_i<q_i.$	The equations form a regular sequence: at each stage the next	equation is monic in a new variable and hence a nonzerodivisor.	Induction on $i$ shows that, for	$Q=\prod_iq_i$, the $Q$-th power of every $Z_i$ lies in $A$.	Thus every element of the quotient has $Q$-th power in $A$.	There is a unique maximal ideal above $(X)$, namely the origin	inverting an element outside it preserves the $Q$-power assertion.	Consequently $R^Q\subseteq A$, and the inclusion $A\to R$	induces an isomorphism on perfections. Nilpotents, if present,	disappear on perfection.	The same finite free presentation over $\widehat A$ proves the	completed assertion.
 \end{proof}

 The following extends \cite[Example 3.5]{Shim} who worked with $n=4$:

 \begin{example}\label{thm:nonCM}	Let $k$ be a perfect field of characteristic $p>0$, and let	$n=p^e\geq4$. Set	$$	S:=k[[u,v]],\qquad
 	a:=u^n,\quad b:=u^{n-1}v,\quad c:=uv^{n-1},\quad d:=v^n,	$$	and let
$	R_n:=k[[a,b,c,d]]\subseteq S.$	Then:	\begin{enumerate}	\item $R_n$ is a complete Noetherian local domain with	$\dim R_n=2$ and $\edim R_n=4$.	\item $R_n^\infty=S^\infty$, so its perfection is coherent.	\item $\depth R_n=1<2=\dim(R_n)$.
 		\item Its normalization is
$B_n=k[[u^n,u^{n-1}v,u^{n-2}v^2,\ldots,uv^{n-1},v^n]], $	and $\edim B_n=n+1$.
 	\end{enumerate}
 	Thus, even with $\edim R_n=4$ and coherent perfection, the
 	embedding dimension of the normalization is unbounded as $e$
 	increases.
 \end{example}
 
 \begin{proof}
 Let  $A:=k[[a,d]]$.	The monomials $b,c$ are integral over $A$, since
$b^n=a^{n-1}d$ and $ c^n=ad^{n-1}.$	The image $k[[a,b,c,d]]=A[b,c]$ is finite over $A$;	a finite algebra over the complete ring $A$ is complete. 	It is a subring of $S$, hence a domain, and has dimension two. 	Its maximal ideal is generated by $a,b,c,d$. 	These four distinct monomials all have ordinary total degree $n$, 	whereas every element of its squared maximal ideal has order at 	least $2n$ in $k[[u,v]]$. Their classes are linearly independent. 	Thus $\edim R_n=4$.
 	 	Since $n$ is a power of $p$,
$	S^n\subseteq k[[u^n,v^n]]=A\subseteq R_n\subseteq S. $	Consequently $R_n^\infty=S^\infty$.
 	Consider the monomial
$	h=u^2v^{n-2}.$	It does not belong to $R_n$: its total degree is $n$, but it is
none of $a,b,c,d$, and the degree-$n$ part of $R_n$ is their
$k$-linear span.
 There is an identity
$	d\,b^{n-2}=a^{n-3}c^2$ and
 $	\frac{b^{n-2}}{a^{n-3}}=h\notin R_n.	$
Hence $b^{n-2}$ represents a nonzero class in
$R_n/a^{n-3}R_n$ annihilated by $d$. 	The pair $a^{n-3},d$ is a system of parameters: $a,d$ already 	is one, by finiteness over $A$. 	It is not a regular sequence, so $R_n$ is not Cohen--Macaulay.	Being a two-dimensional local domain, it has depth at least one; 	therefore its depth equals one.	It remains to verify the normalization assertion. 	The ring $B_n$ consists exactly of those series in $S$ whose 	monomials have total degree divisible by $n$. It is free over $A$ 	with basis  $\{1, t_i=u^{n-i}v^i:1\leq i\leq n-1 \}.	$ 	For $i\geq1$, $t_i=b^i/a^{i-1}$, so the previous formula shows	$B_n\subseteq\operatorname{Frac}(R_n)$. 	It is finite over $R_n$, and their fraction fields agree. 	We prove that $B_n$ is normal directly. 	Decompose a series in $S$ into its $n$ components according 	to total degree modulo $n$. 	If $z\in S\cap\operatorname{Frac}(B_n)$, write $z=f/g$	with $f,g\in B_n$, $g\neq0$. 	The equation $gz=f$, compared component by component, forces	all components of $z$ of nonzero degree modulo $n$ to vanish.	Thus
 $	S\cap\operatorname{Frac}(B_n)=B_n.	$	An element of $\operatorname{Frac}(B_n)$ integral over $B_n$	is integral over $S$ as well, and hence lies in $S$, since	$S$ is normal. It therefore belongs to $B_n$. 	This proves normality and identifies $B_n$ with the normalization.
Finally, the $n+1$ monomials of total degree $n$ generate the	maximal ideal of $B_n$, and are independent modulo its square. 	Thus $\edim B_n=n+1$.
 \end{proof}

 \begin{example}\label{ex:three}	Let $k$ be perfect and
 $R:=\frac{k[[x,y,z,w]]}{(z^p-x^2-y^3-w^5)}.$	Then $\dim R=3$, $\edim R=4$, and
$	R^\infty=k[[x,y,w]]^\infty.$	Indeed the defining equation is monic in $z$, and	$R^p\subseteq k[[x,y,w]]$. 	The equation has no linear term, which gives the asserted embedding 	dimension. It is a domain: $x^2+y^3+w^5$ is not a $p$-th power 	in the fraction field of $k[[x,y,w]]$, since this regular ring is 	normal and the series has a monomial exponent not divisible by $p$. 	Thus $T^p-(x^2+y^3+w^5)$ is irreducible over that field. 	An explicit regular extension is obtained inside $k[[u,v,t]]$ by 	$$	x=u^p,\quad y=v^p,\quad w=t^p,\quad z=u^2+v^3+t^5.$$ 	In particular, this non-binomial hypersurface has coherent perfection	and a regular model.
 \end{example}
 
 \begin{example}
\label{ex:two}
 Let
 $R:=\frac{k[[x,y,z,w]]}	{(z^p-x^{p+1},\ w^p-y^{p+1})}.$	This is a domain of dimension two and embedding dimension four,	and is a complete intersection.	Under	$$x=u^p,\quad y=v^p,\quad z=u^{p+1},\quad w=v^{p+1},	$$	it identifies with	$	k[[u^p,v^p,u^{p+1},v^{p+1}]]\subseteq S=k[[u,v]].	$	To verify injectivity, the quotient displayed  in R is free over	$k[[x,y]]$ on $z^iw^j$, $0\leq i,j<p$; their images have	distinct pairs of exponents modulo $p$, so are linearly independent. 	Since $S^p\subseteq R$, the perfections agree. 	Here $S$ is also the normalization, because	$u=z/x$ and $v=w/y$ lie in the fraction field of $R$.
 \end{example}
 
In a reduced Noetherian ring $D$, write
$ D^\circ=D\setminus\bigcup_{P\in\Min(D)}P.$
The Frobenius closure and tight closure are given by
\begin{align*}	I^F&:=\{u:\ u^q\in I^{[q]}\text{ for some }q=p^e\},\\
	I^*&:=\{u:\ \text{some }c\in D^\circ\text{ satisfies }
	cu^q\in I^{[q]}\text{ for all }q\gg0\}.
\end{align*}

 \begin{example}\label{prop:quartic}	Let $k$ be perfect of characteristic $p\equiv1\pmod4$. Then
$R:=\frac{k[[x,y,z,w]]}{(x^4+y^4+z^4+w^4)}$	is a three-dimensional normal local hypersurface of embedding	dimension four, but $R^\infty$ is not coherent.
 \end{example}
 
 \begin{proof}	Write $f=x^4+y^4+z^4+w^4$ and $r=(p-1)/4$.	The coefficient of $(xyzw)^{p-1}$ in $f^{p-1}$ is
$\frac{(p-1)!}{(r!)^4}\neq0\quad\text{in }k.$	Therefore $f^{p-1}\notin(x^p,y^p,z^p,w^p)$, and Fedder's 	criterion \cite{Fedder} shows that $R$ is $F$-pure. 	The Jacobian ideal has radical the maximal ideal, since $p\neq2$. 	Thus $T$ is regular off its closed point; as a hypersurface it	satisfies $S_2$. Serre's criterion makes it a normal local domain.
Let $I=(y,z,w)R$. The class $x^3$ is not in $I$, since 	$R/I=k[[x]]/(x^4)$. 	We nevertheless have $x^3\in I^*$. For every $q=p^e=4m+1$, the defining equation gives
$x^{3q}=(-1)^{3m}x^3(y^4+z^4+w^4)^{3m}. $ 	For a term with indices $i+j+k=3m$, either one of 	$i,j,k$ is at least $m+1$, or $i=j=k=m$. 	In the first case the term already belongs to 	$(y^q,z^q,w^q)R$. 	In the second case multiplication by $y$ puts it there. 	Hence the fixed nonzero multiplier $y$ satisfies	$y(x^3)^q\in I^{[q]}$ for every  $q=p^e $ 	which proves $x^3\in I^*$. 	Since $R$ is $F$-pure, $I^F=I$, so $I^*\neq I^F$. 	But coherence of the perfection forces equality of tight closure 	and Frobenius closure \cite[Section~3]{Shim}. 	This contradiction proves that $R^\infty$ is not coherent.
 \end{proof}
 
\section{Dealing with(out) excellent-like assumption}
 
 We start by removing an excellent assumption from \cite{Shim}.
 \begin{proposition}\label{prop:one}	Let $R$ be a one-dimensional Noetherian local domain of characteristic	$p>0$. If $A=R^\infty$ is coherent, then the normalization $B$ of $R$ in	$\Frac(R)$ is a discrete valuation ring and
$R\subseteq B\subseteq A$ with $ B^\infty=A.$
 \end{proposition}
 
 \begin{proof}	A coherent perfection of a Noetherian local domain is a GCD domain,	and hence integrally closed; see \cite[Corollary~3.7]{moh}. 	For $b\in B$, the monic equation showing that $b$ is integral over $R$ 	also shows that it is integral over $A$. Since 	$b\in\Frac(R)\subseteq\Frac(A)$, integral closedness gives $b\in A$.	Thus $R\subseteq B\subseteq A$.	Krull--Akizuki implies that $B$ is Noetherian; being the normalization	of a one-dimensional Noetherian domain in its own fraction field,	it is a Dedekind domain.	We must also check that it is local. The ring $A$ is integral over $B$:	every element of $A$ has a $p$-power in $R$. Hence every maximal ideal	of $B$ is the contraction of a maximal ideal of $A$. Perfection of a	local ring is local, so $A$ has just one maximal ideal. Consequently	$B$ has just one maximal ideal and is a DVR.
 Finally, inside $\Frac(R)^\infty$, perfection of the inclusions gives	$A=R^\infty\subseteq B^\infty\subseteq A^\infty=A.$	This proves the claimed equality.
 \end{proof}
 
 The normalization need not be module-finite over $R$ in this argument.
 Its being a Noetherian regular local ring is sufficient for the
 question. The complete one-dimensional case is also treated in
 \cite[Corollary~3.14]{moh}.
 We continue by removing an excellent assumption from \cite{Shim}.
 \begin{proposition}\label{prop:bcm}	Let $(R,\m)$ be a Noetherian local domain of characteristic $p>0$. 	If $R^\infty\cong S^\infty$ for a regular Noetherian local ring $S$, then	$R^\infty$ is a balanced big Cohen--Macaulay $R$-algebra. Consequently,	if every ideal of $R$ is Frobenius closed, then $R$ is Cohen--Macaulay.
 \end{proposition}
 
 \begin{proof}	Identify $S$ with its image in $A=R^\infty$ under the isomorphism.	Then $A=S^\infty$ as an extension of this copy of $S$.	The homeomorphisms on spectra give	$\dim R=\dim S=d$.	Let $x_1,\ldots,x_d$ be any system of parameters of $R$.	Choose a common $q=p^e$ such that all $x_i^q$ belong to $S$.	Their ideal in $S$ is primary for its maximal ideal: a prime of $S$	containing them lifts to a prime of $A$ containing every $x_i$, and	its contraction to $R$ must be $\m$. Uniqueness of primes under	perfection then forces it to be the maximal prime.
 Thus $x_1^q,\ldots,x_d^q$ is a system of parameters of the regular	local ring $S$, hence an $S$-regular sequence. By Kunz's theorem the 	successive Frobenius root extensions of $S$ are flat, so $S\to A$	is flat. The sequence remains regular on $A$, with proper generated 	ideal because $A$ is local. Applying the inverse of the Frobenius 	automorphism $a\mapsto a^q$ of $A$ shows that	$x_1,\ldots,x_d$ is $A$-regular. This proves the first assertion.	For the second, write $I^F$ for the Frobenius closure of $I$.	The equality
 $IA\cap R=I^F$	follows by raising a finite expression in $IA$ to a common	$p$-power. Under the hypothesis, $IA\cap R=I$ for every ideal $I$.	If $a\in R$ and	$ax_i\in (x_1,\ldots,x_{i-1})R$, regularity on $A$ gives	$a\in(x_1,\ldots,x_{i-1})A$; contraction gives	$a\in(x_1,\ldots,x_{i-1})R$.	Thus the chosen system of parameters is $R$-regular and $R$ is	Cohen--Macaulay.
 \end{proof}
 
 \begin{corollary}\label{cor:weak}	A weakly $F$-regular Noetherian local domain which is not	Cohen--Macaulay has no regular Noetherian local model for its	perfection. The same conclusion holds for an $F$-pure	Noetherian local domain which is not Cohen--Macaulay.
 \end{corollary}
 
 \begin{proof}	In the weakly $F$-regular case,	$I\subseteq I^F\subseteq I^*=I$ for every ideal $I$.	In the $F$-pure case, every ideal is likewise Frobenius closed.	Apply Proposition~\ref{prop:bcm}.
 \end{proof}

 Let
$T=\widehat R, \fn=\fm T,
 A=R^\infty$ and $ B=T^\infty.$
 For an ideal $I$ in a characteristic-$p$ ring, $I^{[q]}$ denotes
 the ideal generated by the $q$th powers of its elements, where
 $q=p^e$. 
 
 \begin{lemma}\label{lem:flat}	If $D\to E$ is a flat homomorphism of characteristic-$p$ rings,	then $D^\infty\to E^\infty$ is flat. In particular, $A\to B$ is	faithfully flat.
 \end{lemma}
 
 \begin{proof}	Consider the two Frobenius direct systems for $D$ and $E$.	At every stage the vertical map is the given flat map $D\to E$.	A finite linear relation in $E^\infty$, with coefficients in	$D^\infty$, is represented at a common finite stage, and the 	relation holds after advancing to a further stage. Flatness there	expresses it as a linear combination of relations over that copy	of $D$. Passing to the limit proves the equational criterion for	flatness. This does not require the Frobenius transition maps	themselves to be flat.
 Apply this to  $R\to T$. The induced map	of perfections is local: an element represented at a finite stage	is a unit precisely when its representative is a unit at that	stage. A flat local homomorphism of nonzero local rings is	faithfully flat.
 \end{proof}
 
 \begin{lemma}\label{lem:colon}	Suppose that $A$ is coherent. If $Q\subset T$ is $\fn$-primary	and $u\in T$, then $(QB:_B u)$ is finitely generated over $B$.
 \end{lemma}
 
 \begin{proof}	Choose $N$ with $\fn^N\subseteq Q$. The isomorphism
 	$R/\fm^N\simeq T/\fn^N$ shows that, for $I=Q\cap R$,
$IT=Q$ and $ R/I\simeq T/Q.$
For completeness, $I$ contains $\fm^N$; after reduction modulo	$\fn^N$, the image of $I$ is exactly $Q/\fn^N$, proving $IT=Q$.	Choose $a\in R$ with $u-a\in Q$. Then	$$	(QB:_B u)=(IB:_B a)=(IA:_A a)B.	$$	The last equality follows from flatness of $A\to B$, by tensoring	the kernel of the multiplication map $A\to A/IA$, $z\mapsto az$.	The ideal $IA$ is finitely generated, so coherence of $A$ implies	that $(IA:_A a)$ is finitely generated. Its extension to $B$ is
therefore finitely generated.
 \end{proof}
 
 \begin{lemma}\label{lem:roots}	Let $(D,\mathfrak d)$ be a reduced Noetherian local ring of	characteristic $p>0$. A proper finitely generated ideal of	$D^\infty$ cannot contain $c^{1/p^e}$ for all sufficiently large
$e$, for a fixed $c\in D^\circ$.
 \end{lemma}
 
 \begin{proof}
 Suppose $J=(b_1,\ldots,b_s)$ is such an ideal. Choose a minimal
 prime $P$ of $D$ and put $E=D/P$. The quotient map induces a	surjective local map $D^\infty\to E^\infty$. Its image $J'$ of	$J$ is proper, and the image $\bar c$ of $c$ is nonzero.	If $E$ is a field, then $J'=0$, already a contradiction.
 Otherwise choose a discrete valuation on $\Frac(E)$ whose
valuation ring dominates $E$; such a valuation exists. Extend it to $\Frac(E)^\infty$ by
$	v(z^{1/p^e})=p^{-e}v(z). $	It is nonnegative on $E^\infty$ and strictly positive on every	nonzero nonunit of $E^\infty$. Among the nonzero images of the	$b_i$, let $\delta>0$ be the minimum valuation. There must be at	least one such image because $\bar c^{1/p^e}\in J'$.	Every nonzero element of $J'$ has valuation at least $\delta$,	whereas
 $v(\bar c^{1/p^e})=p^{-e}v(\bar c)\longrightarrow0.$ This contradiction proves the lemma.
 \end{proof}
 
 \begin{proposition}\label{thm:completion}	Assume that $T=\widehat R$ is reduced and $R^\infty$ is coherent.	Then	$$	Q^*=Q^F	\qquad\text{for every $\fn$-primary ideal $Q\subset T$}.	$$	Consequently, if every ideal of $T$ is Frobenius closed, then $T$	is weakly $F$-regular.
 \end{proposition}
 
 \begin{proof}	Suppose $u\in Q^*\setminus Q^F$. Choose $c\in T^\circ$ such that	$cu^q\in Q^{[q]}$ for every sufficiently large $q=p^e$.	The ideal
 $J=(QB:_B u)$ is proper, since $QB\cap T=Q^F$, and is finitely generated by	Lemma~\ref{lem:colon}. Taking $q$th roots in $B$ gives
$c^{1/q}u\in QB$ and $ c^{1/q}\in J$ for $ q\gg0$.	Lemma~\ref{lem:roots} gives a contradiction. Thus $Q^*\subseteq	Q^F$, and the reverse inclusion always holds.
 If every ideal of $T$ is Frobenius closed, all $\fn$-primary	ideals of $T$ are therefore tightly closed. For an arbitrary	proper ideal $L\subset T$ and $u\in L^*$, monotonicity gives	$ u\in (L+\fn^N)^*=L+\fn^N$ { for every } 	$N\geq1.$	Krull intersection in $T/L$ yields $u\in L$. Hence every ideal
 of $T$ is tightly closed.
 \end{proof}
 
 \begin{remark}	The root-and-valuation argument is the same method used in	Shimomoto's tight-closure criterion	\cite[Proposition 3.12; Proposition 3.11]{Shim}.	The extra point here is Lemma~\ref{lem:colon}. 
 \end{remark}
 
 \begin{lemma}\label{lem:fclosed}	If every ideal of $R$ is Frobenius closed, the same holds for $T$. 	In particular, $T$ is reduced.
 \end{lemma}
 
 \begin{proof} 	First let $Q\subset T$ be $\fn$-primary and write $Q=IT$, with 	$I=Q\cap R$, as in Lemma~\ref{lem:colon}. If $u\in Q^F$, choose 	$a\in R$ with $u-a\in Q$. For some $q=p^e$,
$a^q\in Q^{[q]}=I^{[q]}T.$	Faithful flatness of $R\to T$ gives $a^q\in I^{[q]}$, hence	$a\in I$ and $u\in Q$. Thus $Q^F=Q$ 	For any proper ideal $L\subset T$ and $u\in L^F$, it follows that	$u\in (L+\fn^N)^F=L+\fn^N$ for every $N$. Krull intersection	gives $u\in L$. Finally, $(0)^F=(0)$ excludes nonzero nilpotents.
 \end{proof}
 The following drops  a reduced assumption \cite{Shim}.
 \begin{corollary}\label{cor:CM} 
 	Let $(R,\fm)$ be a Noetherian local ring of characteristic $p>0$. 	If $R^\infty$ is coherent and every ideal of $R$ is Frobenius
 	closed, then $\widehat R$ is weakly $F$-regular and $R$ is
 	Cohen--Macaulay. In particular,
 	$$	\left.
 	\begin{array}{l}
 		R\text{ weakly $F$-regular},\\
 		R\text{ not Cohen--Macaulay}
 	\end{array}\right\}
 	\quad\Longrightarrow\quad R^\infty\text{ not coherent}.$$
 \end{corollary}
 
 \begin{proof}	Combine Lemma~\ref{lem:fclosed} with	Theorem~\ref{thm:completion}. A complete weakly $F$-regular local	ring is Cohen--Macaulay: it is a quotient of a regular local ring,	so \cite[Proposition 4.2(c)]{HH} applies.	Cohen--Macaulayness descends from $\widehat R$	to $R$. If $R$ is weakly $F$-regular, then	$I\subseteq I^F\subseteq I^*=I$ for every ideal $I$, so the	hypothesis on Frobenius closure is automatic.
 \end{proof}

 \begin{example}\label{thm:explicit}
Let $R$ be so that $ T=\widehat R\simeq k[[x,y,z,w]]/(xy,xz)
$ and elements $a,b\in R$ whose images in $T$ are
$a=u(x+w), b=vw,
  u,v\in T^\times.$
The ideal
 $	H=(aR^\infty:_{R^\infty}b)$	is not finitely generated. In particular, $R^\infty$ is neither	coherent nor a GCD domain.
 \end{example}
 
 \begin{proof}	Set $A=R^\infty$ and $B=T^\infty$. Lemma~\ref{lem:flat} and	flat base change for a principal colon give 	$	HB=(aB:_B b)=((x+w)B:_B w)=:J.$	We show directly that $J$ is not finitely generated.	First, $J$ is proper. The substitution $y=z=0$, $w=-x$ induces 	a map $B\to k[[x]]^\infty$ which kills $x+w$ and sends $w$ to 	the nonzero element $-x$. Hence $w\notin(x+w)B$.	For every $q=p^e$, one has $xy^{1/q}=0$ in $B$: its $q$th power	is $x^qy=0$, and $B$ is reduced. Consequently,	$wy^{1/q}=(x+w)y^{1/q},$ 	with $y^{1/q}\in J.	$	Now, use the surjective local map	$\pi:B\longrightarrow C=k[[y]]^\infty$ via $ x,z,w\longmapsto0.$	The ideal $\pi(J)$ is contained in the maximal ideal $\mathfrak c$	of $C$, since $J$ is proper and $B$ is local. On the other hand,	it contains $y^{1/p^e}$ for every $e\geq0$. These elements	generate $\mathfrak c$: every nonunit of $C$ is $f^{1/q}$ for	some $f\in yk[[y]]$, and is therefore divisible by $y^{1/q}$.	Thus
 $\pi(J)=\mathfrak c.$	The nonzero ideal $\mathfrak c$ is not finitely generated.	Indeed, perfection implies	$\mathfrak c=\mathfrak c^2$: if $h\in\mathfrak c$, then	$h=(h^{1/p})^p\in\mathfrak c^2$. If $\mathfrak c$ were finitely	generated, Nakayama's lemma would give $\mathfrak c=0$, a	contradiction. Therefore $J$ cannot be finitely generated.
 If $H$ were finitely generated over $A$, its extension $J=HB$	would be finitely generated over $B$. Hence $H$ is not finitely	generated. Coherent rings have finitely generated colons of	finitely generated ideals by single elements. Moreover, in a	GCD domain a colon of two nonzero principal elements is principal.	Both conclusions follow.
 \end{proof}

\section{F-coherent ring without regular model}
Throughout this section $k:=\F_p$. We work with
$$\begin{aligned}
\mathcal{O}&:=k[[u,v]], & T&:=k[[u^2,u^3,v]]=k[[v]]+u^2	\mathcal{O},
&P:=u^2	\mathcal{O}\subset T.
\end{aligned}
$$ 
 The ring $T$ is complete, with maximal ideal
$\n=(u^2,u^3,v)T$. It is isomorphic to
$k[[X,Y,V]]/(Y^2-X^3)$ in every characteristic.
We first observe that
$T/P\cong k[[v]].$
In particular, $P$ is a prime ideal of $T$. Also,
\(
\operatorname{Frac}(T)=\operatorname{Frac}(\mathcal O)=H,
\)
since $u^2\mathcal O\subseteq T$ and
$u=\frac{u^3}{u^2}\in\operatorname{Frac}(T).$
The element $u$ is integral over $T$, and we have:
\begin{fact}\label{eq:cusp}
$\mathcal{O}\cap T_P=T$ and $ u\notin T_P$.
 	\end{fact}

\begin{proof}
To see this, first assume $g\in T\setminus P$ and $h\in 	\mathcal{O}$ with $gh\in T$.	Write an arbitrary element $g\in T\setminus P$ as
$g=a_0(v)+u^2g_2,$	where $a_0(v)\in k[[v]]\setminus\{0\}$ and $g_2\in\mathcal O$.	Similarly, write
$h=b_0(v)+b_1(v)u+u^2h_2,$	where $b_0(v),b_1(v)\in k[[v]]$ and $h_2\in\mathcal O$.	Since $gh\in T=k[[v]]+u^2\mathcal O$, its coefficient	of $u$ must be zero. On the other hand, this coefficient is	$a_0(v)b_1(v)$. Thus
$a_0(v)b_1(v)=0.$	Since $k[[v]]$ is a domain and $a_0(v)\neq 0$, we obtain	$b_1(v)=0$, and hence $h\in T$.
Now, if $h\in\mathcal O\cap T_P$, then $h=a/g$ for some	$a\in T$ and $g\in T\setminus P$. In particular, $gh=a\in T$,	so the preceding argument shows that $h\in T$. Therefore
$\mathcal O\cap T_P=T.$	Finally, since $u\in\mathcal O$ but $u\notin T$, the equality	above implies that $u\notin T_P$.
\end{proof}

\begin{lemma}\label{lem:model}
Let $R\subseteq T$ be a Noetherian local subring with residue field
 $k$ whose inclusion induces $\widehat R=T$. If
$R\perf\simeq S\perf$ for some regular Noetherian local ring $S$,
then there are $b_1,b_2\in R$ and $Q=p^e$ for which
$	T[b_1^{1/Q},b_2^{1/Q}]$
is a complete regular local ring inside $H\perf$.
\end{lemma}
\begin{proof} Transport $S$ into $R\perf$ by the isomorphism. We then have
$S\subseteq R\perf=S\perf\subseteq T\perf=	\mathcal{O}\perf$, with the last
equality following from $	\mathcal{O}^p\subseteq T\subseteq 	\mathcal{O}$.
All the inclusions in this display are local. Perfection induces a
homeomorphism of prime spectra. Hence $S$ has dimension two, and
its residue field is $k$, since it is a subfield of $k=\F_p$.
	Choose regular parameters $s_1,s_2$ of $S$ and a power $q$ of
$p$ with $b_i:=s_i^q\in R$. The only prime of $R$ containing
both $b_i$ is the maximal ideal, by the same homeomorphisms.
Thus $b_1,b_2$ are parameters in $R$ and in $T$.
Evaluation identifies $B=k[[b_1,b_2]]$ with a power series subring
of $T$, and $T$ is finite over $B$. More explicitly,
$R\to T/(b_1,b_2)T$ is onto by completion. Choose
$r_1,\ldots,r_h\in R$ lifting a $k$-basis of this quotient.
Successive expansion modulo powers of $(b_1,b_2)T$ expresses every
element of $T$ as a linear combination of the $r_j$ over the
image of $k[[Y_1,Y_2]]$. Completeness justifies these expansions.
The kernel of evaluation is zero by finiteness and $\dim T=2$.
Inside $T\perf$, the ring $B^{1/q}$ is complete regular with
parameters $s_1,s_2$ and coefficient field $k$. We claim that
$S\subseteq B^{1/q}$ in this same ambient ring. For $s\in S$,
expand $s$ in $\widehat S=k[[s_1,s_2]]$ and let $f_N$ be its
polynomial truncation in total degrees $<N$. Evaluate the whole
series in $B^{1/q}$, obtaining $s'$. Both $s-f_N$ and
$s'-f_N$ belong to $(s_1,s_2)^N 	\mathcal{O}\perf$. The order valuation
of $	\mathcal{O}$ extends to $\mathcal{O}\perf$ by
$\ord(z^{1/p^a})=p^{-a}\ord(z)$. It is positive on nonzero
nonunits. Every nonzero element of $(s_1,s_2)^N 	\mathcal{O}\perf$ has order
at least $N\min_i\ord(s_i)>0$. Taking all $N$ gives $s=s'$.
It follows that $R\subseteq S\perf\subseteq B\perf$. Choose a common
power $Q$ of $p$ with $r_j^Q\in B$ for every $j$. As the
$r_j$ generate $T$ over $B$, this gives $T^Q\subseteq B$.
In particular $B\subseteq T\subseteq B^{1/Q}$. Grouping monomials
of a formal series according to their two exponents modulo $Q$ gives
$B^{1/Q}=B[b_1^{1/Q},b_2^{1/Q}]
 		=T[b_1^{1/Q},b_2^{1/Q}].$
The left side is a complete regular local ring.
\end{proof}

Let $D_1:=\partial/\partial u$ and $D_2:=\partial/\partial v$ be
the formal partial derivatives on $\mathcal{O}$, extended to $H$.
They are ordinary $k$-derivations, so they act on absolute
differentials over $k$.
 	
\begin{lemma}\label{lem:fields}
Suppose $\Aset\subset H$ is a family for which all the elements in
 the indexed family
$(a,D_1a,D_2a)_{a\in\Aset}$
 are algebraically independent over $k$, and let $K:=k(\Aset)$.
Then, for every finite purely inseparable extension $L/K$ generated by
at most two elements of $K\perf$, the multiplication map
$	H\otimes_K L\to HL$
is an isomorphism of fields.
 \end{lemma}
 	
\begin{proof}
Choose $a_1,\ldots,a_n\in\Aset$ so that a sufficiently high
$p$-power of each of the two generators belongs to
$E=k(a_1,\ldots,a_n)$. Let $M/E$ be the finite purely inseparable
 extension generated by the same elements. Write $N=HM$ in $H\perf$.
We first prove that $H\otimes_E M=N$.
Consider the $M$-linear map
$	M^n\rightarrow\Omega_{M/k}$ via $ e_i\mapsto da_i$
and denote its kernel by $W$. Since $k$ is perfect and $M/k$
is finitely generated of transcendence degree $n$, we have
$\dim_M\Omega_{M/k}=n$. The cokernel is $\Omega_{M/E}$, which
has dimension at most two, as $M$ is generated by at most two
elements over $E$. Consequently $\dim_M W\leq 2$.
Let $J=(D_j a_i)_{1\leq j\leq 2,\,1\leq i\leq n}$. Its entries
are algebraically independent over $E$ by $	(a,D_1a,D_2a)_{a\in\Aset}$,
and hence over $M$, which is algebraic over $E$.
Thus $J$ is injective on $W\otimes_M N$. Indeed, if $s=\dim W>0$,	a basis matrix for $W$ has $s$ independent rows. A coordinate
projection to these rows has rank $s$ on $W$. The corresponding
 rank condition for a matrix of indeterminates has a nonzero minor
 polynomial over $M$, which stays nonzero on substituting $J$.
 Let $D:=H\otimes_k M$ and let $\mathfrak p$ be the kernel of
 multiplication $f:D\to N$.
 We identify \(D_{\mathfrak p}\) as the local ring of a smooth
 \(n\)-dimensional variety at a closed point, so that it is regular of
 	dimension \(n\). This is what supplies
\(\dim\mathfrak m_{\mathfrak p}/\mathfrak m_{\mathfrak p}^2=n\).
\(H=k(u,v)\) is the fraction field of \(	\mathcal{O}\), a finitely
generated field of transcendence degree \(2\) over \(k\). \(M\) is a
finite purely inseparable extension of \(E=k(a_1,\dots,a_n)\), hence a
finitely generated field of transcendence degree \(n\) over \(k\).
The image of \(f\) is  a domain, since it is generated over the field
	\(H\) by the field \(M\), so \(HM\) is a field and hence a domain.
Therefore \(\mathfrak p\) is prime and
\(D/\mathfrak p\cong N\) is the residue field of \(D_{\mathfrak p}\).
Since \(M/k\) is finitely generated of transcendence degree \(n\), it is
the function field of an integral \(k\)-variety of dimension \(n\).
Shrinking to an affine open model (using that \(k\) is
perfect), we may write \(M=\Frac(A)\) with \(A\) a smooth integral
\(k\)-algebra of dimension \(n\). Smoothness gives that
\(\Omega_{A/k}\) is locally free of rank \(n\), hence
\(\dim_M\Omega_{M/k}=n\).
Now \(D=H\otimes_k M=H\otimes_k\Frac(A)\). The prime \(\mathfrak p\)
corresponds to the generic point of the base-changed variety
\(X_H=\Spec(H\otimes_k A)\): the map \(D\to N=HM\) is the map to the
function field of \(X_H\), so \(\mathfrak p\) is the zero ideal of
\(H\otimes_k A\) after localizing \(A\). The localization
\(D_{\mathfrak p}\) is the local ring of \(X_H\) at that point.
Because \(N=HM\) is a \emph{finite} field extension of \(H\) (as
\(M/E\) is finite and \(E\subseteq H\)), the quotient \(D/\mathfrak p=N\)
is a finite \(H\)-module, so \(\mathfrak p\) is a \emph{closed} point of
	\(\Spec D\). Thus \(D_{\mathfrak p}\) is the local ring at a closed
point of the \(H\)-variety \(X_H\).
Finally, \(X_H\) is the base change of the smooth
\(k\)-variety \(\Spec A\) along \(\Spec H\to\Spec k\). Since \(k\) is
perfect, \(H/k\) is separably generated, so base change preserves
 smoothness: \(H\otimes_k A\) is smooth over \(H\) of relative dimension
\(n\). A localization of a smooth \(H\)-algebra at a closed point is a
regular local ring of dimension \(n\). Hence \(D_{\mathfrak p}\) is
 		regular of dimension \(n\).
The elements $\delta_i:=a_i\otimes 1-1\otimes a_i$ lie in
 $\mathfrak p$. Their differentials are independent over $N$ in
 $	\Omega_{D/k}\otimes_D N
\simeq(\Omega_{H/k}\otimes_H N)\oplus(\Omega_{M/k}\otimes_M N).	$ 
In a relation among them, the second component puts the coefficient
vector in $W\otimes_M N$. Applying $D_1,D_2$ to the first
component puts it in $\ker J$, so it is zero. The classes of the
$\delta_i$ in the cotangent space of $D_{\mathfrak p}$ are
therefore independent as well, since the differential of a product
of two elements of $\mathfrak p$ vanishes after reduction to $N$.
There are $n$ of them.  		The cotangent space of \(D_{\mathfrak p}\) is
$\mathfrak pD_{\mathfrak p}/\mathfrak p^2D_{\mathfrak p}
 		\cong \Omega_{D_{\mathfrak p}/k}\otimes_{D_{\mathfrak p}}N .
$ The images of the \(\delta_i\) in this space are the classes of
 \(d\delta_i\), and the proof establishes that these are
 \(N\)-linearly independent. Because \(D_{\mathfrak p}\) is regular of
dimension \(n\), its cotangent space has dimension exactly \(n\).
Hence the \(n\) independent classes of the \(\delta_i\) form a
\emph{basis} of the cotangent space.
Let \(\mathfrak m_{\mathfrak p}\) be the maximal ideal of the local ring
	\(D_{\mathfrak p}\). The cotangent space is
\(\mathfrak m_{\mathfrak p}/\mathfrak m_{\mathfrak p}^2\), and we have
just shown that the images of \(\delta_1,\dots,\delta_n\) generate it.
Since \(\delta_i\in\mathfrak m_{\mathfrak p}\), this says
$\mathfrak m_{\mathfrak p}=\mathfrak pD_{\mathfrak p}
=(\delta_1,\dots,\delta_n)+\mathfrak m_{\mathfrak p}^2 .$
Equivalently,
$\mathfrak pD_{\mathfrak p}
=\mathfrak pD_{\mathfrak p}^2+(\delta_1,\dots,\delta_n).$
Nakayama's lemma applied to the finitely generated module
\(\mathfrak pD_{\mathfrak p}\) over the local ring \(D_{\mathfrak p}\)
 	yields
$\mathfrak pD_{\mathfrak p}=(\delta_1,\ldots,\delta_n)D_{\mathfrak p}$.
Now $D/(\delta_1,\ldots,\delta_n)=H\otimes_E M$, since the $a_i$
generate $E$ as a field. This finite $H$-algebra has just one
prime, by pure inseparability of $M/E$. Its localization at that
prime is the residue field just computed. Hence it is $N$.
Finally $K\otimes_E M$ embeds in $H\otimes_E M$. It is a finite
	domain over $K$, hence the field $L=KM$. Tensoring with $H$
over $K$ proves the assertion.
 	\end{proof}

\begin{proposition}\label{prop:ring}
	There exists a countable family $\mathcal A\subseteq T$ such that
$\bigl(a,D_1a,D_2a\bigr)_{a\in\mathcal A}
	$
	is algebraically independent over $k$, 
	and, writing
$	K:=k(\mathcal A)$ and  $R:=T\cap K,$
	we have
	\[
	K\subseteq T_P,\qquad
	\operatorname{Frac}(R)=K,\qquad
	\widehat R\cong T.
	\]
	Moreover, $R$ is a countable normal Noetherian local domain,
	and
$P\cap R=0.$
	In particular,
$\dim R=2$ and $
	\operatorname{edim}R=3.$
\end{proposition}

\begin{proof}
	We divide the construction into several steps.

	\text{Step 1. The field-enlargement construction.}
For $z\in T$, write
$\overline z=z(0,v)\in k[[v]].$
We construct a countable family $\mathcal A\subseteq T$ satisfying
the following two conditions:
\begin{enumerate}
\item The indexed family
$(a,D_1a,D_2a)_{a\in\mathcal A}	$
is algebraically independent over $k$.
\item The family of residues
$(\overline a)_{a\in\mathcal A}$
is algebraically independent over $k$.
\end{enumerate}
	
The second condition will ensure that the field generated by
$\mathcal A$ is contained in $T_P$.
The key observation is the following. Suppose that
$\mathcal A_0\subseteq T$ is a countable family satisfying both
conditions. Let $F\subseteq H$ be a countable subfield containing
$	k,\ u,\ v,\ b,\ c,$
the elements of $\mathcal A_0$, their two derivatives, and
their residues, where
$	c\in T$ and $ b\in T\setminus P.$
We can choose
$t_0,t_1,t_2\in k[[v]]$
algebraically independent over $F$. Indeed, $F$ is countable,	whereas $k[[v]]$ is uncountable, and the algebraic closure of
a countable field is countable.
Define
$	h:=t_0^p+u^{p+1}t_1^p+vt_2^p$	and  $ a=c+bh$.
	Since $h\in T$, we have
$	a\in c+bT.$ In characteristic $p$,
$	D_1h=u^pt_1^p$ and $D_2h=t_2^p.$
Consequently,	\[	\begin{aligned}		a&=c+bh,\\		D_1a&=D_1c+(D_1b)h+bu^pt_1^p,\\		D_2a&=D_2c+(D_2b)h+bt_2^p.\end{aligned}	\]	Thus the transformation from	\[	(t_0^p,t_1^p,t_2^p)	\quad\text{to}\quad	(a,D_1a,D_2a)	\]	is affine linear over $F$, with determinant
$	b^3u^p\neq0.$	It follows that the new triple is algebraically independent
over the field generated by the previous triples. Hence adjoining	$a$ preserves the first independence condition.
Moreover,	$	\overline a	=\overline c+\overline b(t_0^p+vt_2^p).	$	Since $b\notin P$, we have $\overline b\neq0$. Thus	$\overline a$ is transcendental over the field generated by	the preceding residues. This preserves the second independence	condition.	We have therefore proved the following useful assertion:
\medskip
\noindent	\emph{Given a countable family satisfying the two independence
we can adjoin an element $a\in c+bT$ while preserving both
	independence conditions.}
This is the field-enlargement mechanism used below. At each
application, the field generated by the family is enlarged by
adjoining  new element $a$.
	
	\medskip
	\noindent
	\text{Step 2. Initial elements and the field $K$.}
Let
$	\mathfrak n=(u^2,u^3,v)T$	be the maximal ideal of $T$.	Starting with the empty family, apply the preceding construction	successively with	\[	(c,b)=(u^2,v^2),\quad	(u^3,v^2),\quad	(v,v^2).	\]	Notice that $v^2\notin P$. We obtain three elements	$b_1,b_2,b_3\in T$ such that	\[	b_1\equiv u^2\pmod{\mathfrak n^2},	\qquad	b_2\equiv u^3\pmod{\mathfrak n^2},	\qquad	b_3\equiv v\pmod{\mathfrak n^2}.	\]	We include these elements in our family.
For any countable family $\mathcal A_0$ satisfying the two	independence conditions, let	$	K_0::=k(\mathcal A_0)$ and $R_0:=T\cap K_0.$	The algebraic independence of the residues implies that every	nonzero polynomial in finitely many elements of $\mathcal A_0$	has nonzero residue modulo $P$. Consequently, every nonzero	element of $K_0$ is a unit in $T_P$. In particular,	$$(\ast):\quad 	K_0\subseteq T_P,\qquad P\cap R_0=0$$
	
	We now enlarge the family further in order to arrange an
	ideal-contraction property.
	
	\medskip
	\noindent	\text{Step 3. Enforcing the ideal-contraction property.}
We claim that we can construct the family so that, for	$	K:=k(\mathcal A)$ and $R:=T\cap K,$	we have	$$(+)\quad		IT\cap  R=I		\quad\text{for every finitely generated ideal }I\subseteq R.	$$
		We explain the recursive construction.
Suppose that $\mathcal A_0$ is a countable family satisfying
	the two independence conditions, and set
$K_0:=k(\mathcal A_0)$ and $ R_0:=T\cap K_0.$
	Consider elements
$	c,y_1,\ldots,y_m\in R_0$	such that
$	c\in(y_1,\ldots,y_m)T.$
	We want to enlarge the field so that this membership already
	holds in the corresponding intersection ring.
Choose $z_1,\ldots,z_m\in T$ with
$	c=\sum_{i=1}^m y_iz_i.$
If all $y_i$ are zero, then $c=0$, and there is nothing to prove.
Otherwise, after reordering, we may assume $y_m\neq0$.
By $(\ast)$, every nonzero element of $R_0$ is a	unit in $T_P$. Hence
$y_m\notin P.$
For each $i<m$, apply the construction of Step 1 with	$	c=z_i$ and $b=y_m.	$
	This gives an element
$z_i'\in z_i+y_mT$
	which we adjoin to the family, preserving both independence
	conditions.
	Now define
$	z_m'=	\frac{c-\sum_{i=1}^{m-1}y_iz_i'}{y_m}.$	Since
$	z_i'-z_i\in y_mT,$	we have	\[	c-\sum_{i=1}^{m-1}y_iz_i'	=	y_m\left(z_m-\sum_{i=1}^{m-1}	\frac{y_i(z_i'-z_i)}{y_m}\right).	\]	Thus $z_m'\in T$. Also, all the elements appearing in its	definition belong to the enlarged field. Therefore $z_m'$	belongs to that field as well.
It follows that, in the enlarged intersection ring,	$c=\sum_{i=1}^m y_iz_i'.$	Thus the given ideal-contraction problem has been solved.
Notice that for a principal ideal, no new elements are needed:	if $c\in y_1T$ and $c,y_1\in K_0$, then	$c/y_1\in T\cap K_0$.
Since $R_0$ is countable, there are only countably many finite
	tuples to consider. We can therefore enumerate all the tuples
$	(c;y_1,\ldots,y_m)$	satisfying the above membership condition and solve them	successively, taking the union of the resulting countable	families. The resulting family is still countable and satisfies	both independence conditions.
We repeat this procedure countably many times. More precisely,	we obtain an increasing sequence	\[	\mathcal A_0\subseteq\mathcal A_1	\subseteq\mathcal A_2\subseteq\cdots	\]	of countable families, each satisfying the two independence	conditions, such that at stage $n+1$ every ideal-contraction	problem involving elements of
$	R_n=T\cap k(\mathcal A_n)$
is solved.
Set
	\[
	\mathcal A:=\bigcup_{n\geq0}\mathcal A_n,
	\qquad
	K:=k(\mathcal A),
	\qquad
	R:=T\cap K.
	\]
The family $\mathcal A$ is countable, and both independence
conditions are preserved because every polynomial relation
involves only finitely many elements.
Every finite tuple of elements of $R$ belongs to some $R_n$,
since every element of $K$ is a rational expression involving
only finitely many elements of $\mathcal A$. Hence every
ideal-contraction problem is eventually solved. This proves
$(+)$.

	\medskip
	\noindent
	\text{Step 4. The basic properties of $R$.}	By construction and the independence of the residues,
$	K\subseteq T_P,  P\cap R=0.$
Since $\mathcal A\subseteq R$, we also have
$	\operatorname{Frac}(R)=K.$	We next show that $R$ is local. Set
$\mathfrak m=R\cap\mathfrak n.$
If $r\in R\setminus\mathfrak m$, then $r$ is a unit in $T$.
Since $r^{-1}\in K$ as well, we have
$r^{-1}\in T\cap K=R.$	Thus $R$ is local with maximal ideal $\mathfrak m$.	Moreover, the inclusion $R\subseteq T$ induces an embedding
$R/\mathfrak m\hookrightarrow T/\mathfrak n=k.$ Since $k\subseteq R$, we obtain
$R/\mathfrak m\cong k.$ We now prove that $R$ is Noetherian. Let $J$ be an arbitrary	ideal of $R$. Since $T$ is Noetherian, the ideal $JT$ is	generated by finitely many elements of $J$. Thus there are
$j_1,\ldots,j_r\in J$ such that
$JT=(j_1,\ldots,j_r)T.$	Put $I=(j_1,\ldots,j_r)R$. For every $x\in J$, we have
$	x\in JT=IT.$	Since $x\in R$, it follows from $(+)$
	that
$	x\in IT\cap R=I.$	Therefore $J=I$, and every ideal of $R$ is finitely generated.	Hence $R$ is Noetherian.
Finally, we prove that $R$ is normal. Since	$K\subseteq T_P$, equation Fact 6.1 gives
$	R=T\cap K=\mathcal O\cap K.$	Suppose that $x\in K$ is integral over $R$. Since	$R\subseteq\mathcal O$, the same monic equation shows that	$x$ is integral over $\mathcal O$. But $\mathcal O=k[[u,v]]$	is integrally closed in $H$. Hence $x\in\mathcal O$.	Consequently, $	x\in\mathcal O\cap K=R.$
Thus $R$ is integrally closed in its fraction field $K$.
Therefore $R$ is a normal Noetherian local domain.
Since $\mathcal A$ is countable and $k=\mathbb F_p$ is finite,
	the field $K=k(\mathcal A)$ is countable. Hence $R$ is countable.
	
	\medskip
	\noindent
	\text{Step 5. Identification of the completion.}
	By our initial choice of $b_1,b_2,b_3\in R$, we have
	\[
	b_1\equiv u^2,\qquad
	b_2\equiv u^3,\qquad
	b_3\equiv v
	\pmod{\mathfrak n^2}.
	\]
	Nakayama's lemma gives
$\mathfrak n=(b_1,b_2,b_3)T.$
	In particular,
$\mathfrak mT=\mathfrak n.$ 	Applying $(+)$ to $\mathfrak m^j$,
we obtain
$	\mathfrak n^j\cap R
	=(\mathfrak m^jT)\cap R
	=\mathfrak m^j	$
	for every $j\geq1$.
	We claim that the natural map
$R\longrightarrow T/\mathfrak n^j$
	is surjective for every $j\geq1$.
Indeed, the residue field of $R$ is $k$, and
	$\mathfrak n=(b_1,b_2,b_3)T$. Hence the monomials of degree
	$i$ in $b_1,b_2,b_3$ span
$
	\mathfrak n^i/\mathfrak n^{i+1}
$
	over $k$. All these monomials belong to $R$. Starting with
	the residue field and successively lifting modulo
$	\mathfrak n,\mathfrak n^2,\ldots,\mathfrak n^j,$ we see that every element of $T/\mathfrak n^j$ has a	representative in $R$. Therefore the natural map
$	R/\mathfrak m^j\longrightarrow T/\mathfrak n^j
$ is an isomorphism: it is surjective by the preceding argument,
	and injective because $	\mathfrak n^j\cap R=\mathfrak m^j.$	Passing to inverse limits gives
$	\widehat R
=\varprojlim_j R/\mathfrak m^j\cong
\varprojlim_j T/\mathfrak n^j=T.$ Finally, since completion preserves dimension and embedding
	dimension for Noetherian local rings, and
$	T\cong k[[X,Y,V]]/(Y^2-X^3),$ we obtain
$	\dim R=\dim T=2$	and
$\operatorname{edim}R
	=\operatorname{edim}T=3.$ This completes the proof.
\end{proof}

	\begin{theorem}\label{thm:main}
	For every prime $p$ there is a countable normal Noetherian local
domain $R$ with
$\widehat R\simeq \F_p[[u^2,u^3,v]]$
such that $R\perf$ is coherent and $R\perf\not\simeq S\perf$ for
every regular Noetherian local ring $S$.
\end{theorem}

\begin{proof} 
Take $R$ from Proposition~\ref{prop:ring}. Dimension and embedding
dimension agree with those of its completion $T$. The inclusions
	$	\mathcal{O}^p\subseteq T\subseteq 	\mathcal{O}$ give $T\perf=	\mathcal{O}\perf$.
The latter is coherent because it is the union of the Noetherian rings
$	\mathcal{O}^{1/p^e}$ with flat transition maps, as is seen by grouping monomials in $u,v$ according to the exponents modulo $p$.
Since $R\to T=\widehat R$ is faithfully flat, $F$-coherence
	descends to $R$ \cite[\S 3]{Shim}. In detail, passage to
perfections preserves flatness by the finite-relation criterion
applied at a common Frobenius stage. The induced map is local, hence
faithfully flat. Kernels of finite matrices over $R\perf$ become
finitely generated after this base change and are therefore already
finitely generated.
Suppose now that a regular model exists. Apply Lemma~\ref{lem:model}
and define	$$U:=T[b_1^{1/Q},b_2^{1/Q}],\qquad
L:=K(b_1^{1/Q},b_2^{1/Q}),\qquad
C:=(R\setminus 0)^{-1}T .	$$On the other hand Lemma~\ref{lem:fields} gives $H\otimes_K L=HL$. From Lemma~\ref{lem:model}, the ring
$	U=T[b_1^{1/Q},b_2^{1/Q}]\subseteq H^\infty$
is a complete regular local ring, and \(U\supseteq T\). In particular
$u\in T\subseteq U.$
Also \(U\subseteq C[b_1^{1/Q},b_2^{1/Q}]\), since \(T\subseteq C\) and the
roots are adjoined on both sides. Localizing at \(R\setminus0\), i.e.
tensoring with \(K=\Frac R\) over \(R\), gives
$(R\setminus0)^{-1}U=C\otimes_K L,$
because \(C=T\otimes_R K\) and \(L=K[b_1^{1/Q},b_2^{1/Q}]\).
To see $(C\otimes_K L)\cap H=C$, choose a \(K\)-basis of \(L\) beginning with \(1\). It stays \(H\)-independent in \(H\otimes_K L=HL\). An element of \(H\) lying in \(C\otimes_K L\) is a \(K\)-linear combination
\(\sum_j c_j\otimes\ell_j\) with \(c_j\in C\subseteq H\); comparing
coefficients in the \(H\)-basis, the coefficients of \(\ell_j\) for
	\(j\neq1\) vanish, and the coefficient of \(\ell_1=1\) lies in \(C\).
Hence
$ 	(C\otimes_K L)\cap H=C .$
Now \(u\in U\), so  
$u\in (R\setminus0)^{-1}U=C\otimes_K L .$
	Since also \(u\in H\), we obtain 
$u\in (C\otimes_K L)\cap H=C$.
But $P\cap R=0$ gives $C\subseteq T_P$,
	contradicting Fact \ref{eq:cusp}.\end{proof}
 	
\begin{remark}
The above example is nonexcellent. 
\end{remark}

\end{document}